\documentclass[12pt]{amsart}

\usepackage{amssymb}
\usepackage{latexsym}
\usepackage{verbatim}
\usepackage{fullpage}
\usepackage{pslatex}

\input {cyracc.def}
\font\ququ=cmr10 scaled \magstep1
\newcommand{\re}[1]{\textrm  (\ref{#1})}

\renewenvironment{proof}
{\noindent {\sl Proof.}\quad }{\hfill
$\square$ \vskip1.1ex\noindent }

\newenvironment{proof*}
{\noindent {\sl Proof.}\quad }{\hfill
$\square$}

\renewcommand{\theequation}{\thesection .\arabic{equation}}
\renewcommand{\thesubsubsection}{\theequation .\arabic{subsubsection}}

\catcode`\@=\active
\catcode`\@=11
\def\@eqnnum{\hbox to
.01pt{}\rlap{\hskip-\displaywidth(\mathbf{\theequation})}}
\catcode`\@=12
\newenvironment{s}[1]
{ \vskip1.2ex \refstepcounter{equation}
\noindent {\bf #1\enspace\theequation .} \begin{sl}}{\end{sl}
\vskip1.1ex\noindent }

\newenvironment{rem}[1]
{ \vskip1.2ex \refstepcounter{equation}
\noindent {\bf #1\enspace \theequation .} }{ \vskip1.1ex\noindent }

\newenvironment{rems}[1]
{ \vskip1.2ex \refstepcounter{subsubsection}
\noindent {\bf \thesubsubsection\enspace {#1}.} }{ \vskip1.1ex\noindent }

\newenvironment{subs}[1]
{\vskip1.2ex \refstepcounter{equation}
\noindent {\bf (\theequation)\quad #1.} }{\quad}

\newcommand {\be}{{\frak b}}
\newcommand {\ce}{{\frak c}}

\newcommand {\g}{{\frak g}}
\newcommand {\h}{{\frak h}}

\newcommand {\te}{{\frak t}}
\newcommand {\ut}{{\frak u}}

\newcommand {\sln}{{\frak sl}_n}

\newcommand {\spp}{{\frak sp}_{2p}}

\newcommand {\son}{{\frak so}_n}
\newcommand {\sop}{{\frak so}_p}
\newcommand {\esi}{\varepsilon}
\newcommand {\ap}{\alpha}

\newcommand {\lb}{\lambda}

\newcommand {\vp}{\varphi}
\newcommand {\ca}{{\mathcal A}}

\newcommand {\ccc}{{\mathcal C}}
\newcommand {\ch}{{\mathcal H}}

\newcommand {\N}{{\mathcal N}}

\newcommand {\cs}{{\mathcal S}}
\newcommand {\VV}{{\Bbb V}}

\newcommand {\cl}{{\mathrm{cl}}}

\newcommand {\hot}{{\mathrm{ht\,}}}

\newcommand {\lev}{{\mathrm{lev\,}}}

\newcommand {\rk}{{\mathrm{rk\,}}}

\newcommand {\GR}[2]{{\textrm{{\bf #1}}}_{#2}}

\newcommand {\ov}{\overline}

\newcommand {\AD}{{\frak Ad}}
\newcommand {\AN}{{\frak An}}
\newcommand {\CL}{{\frak Clus}}

\newcommand {\vno}[1]{\vskip#1 ex\noindent}

\newcommand {\qus}{\hfill $\square$ \vno{1.1}}

\newcommand {\beq}{\begin{equation}}
\newcommand {\eeq}{\end{equation}}

\newfam\Bbbfam\newfam\eufam%
\font\Bbbfont=msbm10 scaled 1200%
\font\olala=msam10 scaled 1200%
\font\frak=eufm10 scaled 1400%
\font\Bbbsmallfont=msbm8%
\textfont\Bbbfam=\Bbbfont\scriptfont\Bbbfam=\Bbbsmallfont%%
\font\euzw=eufm10 scaled 1200%
\font\euac=eufm7 scaled 1200%
\font\euacc=eufm7 scaled 1000%
\textfont\eufam=\euzw\scriptfont\eufam=\euac%
\scriptscriptfont\eufam=\euacc%
\def\frak{\fam\eufam}%
\def\Bbb{\fam\Bbbfam}%

\def\varnothing{\hbox {\Bbbfont\char'077}}
\def\square{\hbox {\olala\char"03}}
\def\bbk{\hbox {\Bbbfont\char'174}}

\def\cyeq{\hbox {\olala\char'064}}

\begin{document}
\setlength{\parskip}{2pt plus 4pt minus 0pt}
\hfill {\scriptsize March 9, 2003}
\vskip1ex
\vskip1ex

%\begin{center}
%{\Large \bf %\parbox[t]{\textwidth}
\title[]{{\sf ad}-nilpotent ideals
of a Borel subalgebra: generators and duality}
%\medskip \\
\author[]{\sc Dmitri I. Panyushev}
\thanks{This research was supported in part by the Alexander von
Humboldt-Stiftung and RFBI Grant no. 01--01--00756}
\maketitle
\begin{center}
{\footnotesize
{\it Independent University of Moscow,
Bol'shoi Vlasevskii per. 11 \\
121002 Moscow, \quad Russia \\ e-mail}: {\tt panyush@mccme.ru }\\
}
\end{center}
%\refstepcounter{section}
In this paper, we develop several combinatorial aspects of the theory
of {\sf ad}-nilpotent ideals.
Let $\be$ be a fixed Borel subalgebra of a simple Lie algebra $\g$.
Following
%%Cellini and Papi
\cite{CP1}, we say that an ideal of $\be$ is
{\sf ad}-{\it nilpotent\/}, if it is contained in $[\be,\be]$.
Let ${\AD}$ of $\AD(\g)$ denote the set of all {\sf ad}-nilpotent ideals
of $\be$. Any $\ce\in\AD$ is completely
determined by the the corresponding
set of roots. More precisely, let $\te$ be a Cartan subalgebra of $\g$
lying in $\be$ and let $\Delta$ be the root system of the pair $(\g,\te)$.
Choose $\Delta^+$, the system of positive roots, such
that the roots of $\be$ are positive.
Then $\ce=\oplus_{\gamma\in I}\g_\gamma$, where $I$ is a suitable subset
of $\Delta^+$ and $\g_\gamma$ is the root space for $\gamma\in\Delta^+$.
In particular, this means that $\AD$ is finite.
%%there are finitely many {\sf ad}-nilpotent $\be$-ideals
%and that any question concerning these ideals can be stated in terms of
%combinatorics of the root system.
Abusing language, we shall say that
such $I\subset \Delta^+$ is an {\sf ad}-nilpotent ideal, too.

In \cite{CP1}, Cellini and Papi proved that there is a bijection between
the {\sf ad}-nilpotent $\be$-ideals and the elements of the
affine Weyl group $\widehat W$ satisfying certain property (see
\re{ideals} below). In our paper, these elements are said to be
{\it admissible\/}.
Using admissible elements, Cellini and Papi established a bijection
between $\AD$ and the points of the coroot lattice lying in a
certain $\rk\g$-dimensional simplex $\tilde D$ with rational vertices
\cite{CP2}.
As a consequence, they
obtained a conceptual proof for the explicit formula giving the
number of {\sf ad}-nilpotent ideals in all simple Lie algebras.

In Section~2, we give a characterization of the generators of
{\sf ad}-nilpotent ideals in terms of
admissible elements (Theorem~\ref{gener}).
%Suppose $I\subset \Delta^+$ is an
%{\sf ad}-nilpotent ideal and $w$ is the corresponding admissible
%element. Then $\gamma\in I$ is a generator if and only if
It is then shown that an ideal $I$ has $k$ generators if and only
if the corresponding lattice point
lies on the face of $\tilde D$ of codimension $k$
(Theorem~\ref{main2}). It is curious
that $\tilde D$ has exactly one integral vertex. We deduce this from
the fact that there is only one {\sf ad}-nilpotent ideal
having $\rk\g$ generators.

In Section~\ref{q-sign}, we consider the `simple root' statistic on $\AD(\g)$,
which assigns to any ideal the number of simple roots in it.
Write $\AD(\g)_i$ for the set of ideals containing exactly $i$ simple roots.
We give recurrent formulas for these numbers and then compute them
for $\GR{A}{p}$ and the exceptional Lie algebras. It is also shown
thet the simple root statistic has the same distribution for
$\GR{B}{p}$ and $\GR{C}{p}$. In case of $\GR{C}{p}$ and $\GR{D}{p}$, we give
conjectural values for $\#\AD(\g)_i$, which are, no doubt, true.
As a consequence of this theory, we observe some similarities
between the {\sf ad}-nilpotent ideals and {\it clusters\/} (see
\cite{cluster} for the latter).
It is shown that the simple root
statistic on $\AD(\g)$ and a certain statistic on the set of
clusters have the same distribution (Theorem~\ref{empir}).

To obtain a closed formula for $\#\AD(\g)_0$
(Proposition~\ref{no-simple}), we exploit a bijection
between the {\sf ad}-nilpotent ideals and the regions of the Catalan
arrangement lying in the dominant chamber, see \cite{shi}.
We show that $I\in\AD(\g)_0$ if and only if the corresponding region is
bounded. In turn, the number of bounded regions of any arrangement
can be counted using a powerful result of Zaslavsky,
once one knows the characteristic polynomial, see Prop.~\ref{zaslavsky}
for details. After this part was written, I learned that the formula
for $\#\AD(\g)_0$ was already obtained, in the same way,
in a recent preprint of Athanasiadis~\cite{ath02}.
The main result of Athanasiadis' preprint is a case-free simple proof
of the formula for the characteristic polynomial of the
Catalan arrangement.

In the last three sections, we consider the statistic that
assigns to an ideal $I\in\AD(\g)$ the number of its generators.
In case of $\g=\sln$, the {\sf ad}-nilpotent ideals are identified with Dyck path
of semilength $n$ and, therefore,
the generating function for this statistic is
the famous {\it Narayana polynomial\/} (of degree $n-1$).
For this reason, we say that the generating function for this statistic for arbitrary $\g$ is a
generalized Narayana polynomial.
Motivated by the fact that the
Narayana polynomial is palindromic, we were searching for a materialization
of this property, i.e., for an involutory mapping (duality)
on $\AD(\sln)$
that takes the ideals with $k$ generators to the ideals with $n-1-k$
generators. For $\sln$, such a materialization does exists, and it has
a number of nice properties, see Section~\ref{sl}.
The nicety of these properties is that their formulation admits immediate
generalization to all simple Lie algebras. We also show that the number of
self-dual ideals in ${\frak sl}_{2m+1}$ equals $C_m$, the $m$-th Catalan
number.
In Section~\ref{classical}, the results concerning duality are
extended to series $B$ and $C$. This clearly implies that the generalized Narayana polynomials for $B$ and $C$ (in fact,
they are equal) are palindromic.
We conjecture that such a duality exists for any
simple Lie algebra.
At least, we show that the generalized Narayana polynomials are
palindromic for the exceptional Lie algebras. General properties of this
conjectural duality are discussed in Section~\ref{vmeste}.
%Considerable part of the theory of {\sf ad}-nilpotent ideals can be
%in the combinatorial context.

{\bf Acknowledgements.} This paper was written during my stay at the
Ruhr-Universit\"at Bochum and
Max-Planck-Institut f\"ur Mathematik (Bonn).
I would like to thank both institutions for hospitality
and excellent working conditions.

\section{Preliminaries on {\sf ad}-nilpotent ideals}
\label{prelim}

\noindent
\begin{subs}{Main notation}
\end{subs}
$\Delta$ is the root system of $(\g,\te)$ and
$W$ is the usual Weyl group. For $\ap\in\Delta$, $\g_\ap$ is the
corresponding root space in $\g$.
% and $e_\ap$ is a nonzero vector in $\g_\ap$.

$\Delta^+$  is the set of positive
roots and $\rho=\frac{1}{2}\sum_{\ap\in\Delta^+}\ap$.

$\Pi=\{\ap_1,\dots,\ap_p\}$ is the set of simple roots in $\Delta^+$.

$\ccc$ \ is the fundamental Weyl chamber.
 \\
 We set $V:=\te_{\Bbb Q}=\oplus_{i=1}^p{\Bbb Q}\ap_i$ and denote by
$(\ ,\ )$ a $W$-invariant inner product on $V$. As usual,
$\mu^\vee=2\mu/(\mu,\mu)$ is the coroot
for $\mu\in \Delta$.

$Q=\oplus _{i=1}^p {\Bbb Z}\ap_i \subset V$ is the root lattice and
$Q^\vee=\oplus _{i=1}^p {\Bbb Z}\ap_i^\vee$  is the coroot lattice.

$Q^+=\{\sum_{i=1}^p n_i\ap_i \mid n_i\in {\Bbb N} \}\subset Q$.
\\
Letting $\widehat V=V\oplus {\Bbb Q}\delta\oplus {\Bbb Q}\lb$, we extend
the inner product $(\ ,\ )$ on $\widehat V$ so that $(\delta,V)=(\lb,V)=
(\delta,\delta)=
(\lb,\lb)=0$ and $(\delta,\lb)=1$.

$\widehat\Delta=\{\Delta+k\delta \mid k\in {\Bbb Z}\}$ is the set of affine
real roots and $\widehat W$ is the  affine Weyl group.
\\
Then $\widehat\Delta^+= \Delta^+ \cup \{ \Delta +k\delta \mid k\ge 1\}$ is
the set of positive
affine roots and $\widehat \Pi=\Pi\cup\{\ap_0\}$ is the corresponding set
of affine simple roots.
Here $\ap_0=\delta-\theta$, where $\theta$ is the highest root
in $\Delta^+$.  The inner product $(\ ,\ )$ on $\widehat V$ is
$\widehat W$-invariant.
\\
For $\ap_i$ ($0\le i\le p$), we let $s_i$ denote the corresponding simple
reflection in $\widehat W$.
If the index of $\ap\in\widehat\Pi$ is not specified, then we merely write
$s_\ap$. %% for the corresponding reflection.
The length function on $\widehat W$ with respect
to  $s_0,s_1,\dots,s_p$ is denoted by $l$.
For any $w\in\widehat W$, we set
\[
\widehat N(w)=\{\ap\in\widehat\Delta^+ \mid w(\ap) \in -\widehat \Delta^+ \} .
\]
%If $w\in W$, then $\widehat N(w)\subset \Delta^+$ and we also write
%$N(w)=\widehat N(w)$ in this case.
Our convention concerning
$\widehat N(w)$ is the same as in \cite{Ko1},\,\cite{long}, but
opposite to that in \cite{CP1},\,\cite{CP2}, so that our
$\widehat N(w)$ is $\widehat N(w^{-1})$ in the sense of Cellini-Papi.
%%\cap \Delta^+$.
%\vskip2ex
\begin{subs}{{\sf ad}-nilpotent ideals}               \label{ideals}
\end{subs}
Throughout the paper, $\be$ is the Borel subalgebra of $\g$ corresponding to
$\Delta^+$
and $\ut=[\be,\be]$. The expression ``{\sf ad}-nilpotent ideal"
or just ``ideal'' always refers to a $\be$-ideal lying in $\ut$.
Let $\ce\subset\be$ be an {\sf ad}-nilpotent ideal. Then
$\ce=\underset{\ap\in I}{\oplus}\g_\ap$
for a subset $I\subset \Delta^+$, which is called the {\it set of roots of\/}
$\ce$.
As our exposition will be mostly combinatorial, an {\sf ad}-nilpotent ideal
will be identified with the respective set of roots.
That is, $I$ is said to be an {\sf ad}-nilpotent ideal, too. Whenever we want to
explicitly indicate the context,  we say that $\ce$ is a {\it geometric\/}
{\sf ad}-nilpotent ideal, while $I$ is a {\it combinatorial\/} {\sf ad}-nilpotent
ideal.
Accordingly, being in combinatorial (resp. geometric) context,
we speak about cardinality (resp. dimension) of an ideal.
In the combinatorial context, the definition of an {\sf ad}-nilpotent ideal
can be stated as follows.
\\
$I\subset\Delta^+$ is an {\sf ad}-nilpotent ideal, if the following  condition
is satisfied: \\
\centerline{
 if $\gamma\in I$, $\nu\in\Delta^+$, and $\gamma+\nu\in\Delta$, then
$\gamma+\nu\in I$. }
%\\[.6ex]
We consider $\Delta^+$ as poset with respect to the standard partial
order `$\cyeq$', i.e., $\mu\cyeq \nu$ if and only if $\nu-\mu\in Q^+$.
Therefore, a combinatorial {\sf ad}-nilpotent ideal
is nothing but a {\it dual order ideal\/} of the poset $(\Delta, \cyeq)$.
An element $w\in \widehat W$ is said to be  {\it admissible},   if it has two properties:

(a) \ $w(\ap)$ is positive for any $\ap\in\Pi$;

(b) \ if $w^{-1}(\ap)$ is negative for an $\ap\in\widehat \Pi$, then
$w^{-1}(\ap)=\gamma-\delta$ for some $\gamma\in\Delta^+$.
\\[.6ex]
By  \cite[Sect.\,2]{CP1}, there is a one-to-one correspondence
between the admissible elements of $\widehat W$ and the
{\sf ad}-nilpotent $\be$-ideals. This
correspondence is obtained as follows:
% Note that our definition of the set
%$\widehat N(w)$ is different form that of \cite{CP1}. Therefore the
%formulas defining the admissible elements are also different.

$\bullet$ \
Given $\ce\in\AD$, consider the members of the descending  central series
$\ce=\ce^1$,
$\ce^k=[\ce^{k-1},\ce]$ ($k\ge 2$)  and the corresponding sets of roots
$I_k$. Clearly, $I_k\supset I_{k+1}$ and $I_m=\varnothing$ for $m\gg 0$.
Set $N_k=\{k\delta-\gamma\mid \gamma\in I_k\}$. Then
$\Phi:= \cup_{k\ge 1}N_k$   is a closed subset of $\widehat \Delta^+$ whose
complement is closed as well, and therefore there is a unique $w\in\widehat W$
such that $\Phi=\widehat N(w)$. This $w$ is the required admissible element.

$\bullet$ \ Conversely, if   $w\in\widehat W$ is admissible, then
$\widehat N(w)=\cup_{k\ge 1}N_k$, where $N_k=\{k\delta-\gamma\mid
\gamma\in I_k\}$
and $I_k\subset\Delta^+$. Then $I_1$ is the set of roots of an
{\sf ad}-nilpotent ideal, say $\ce$.   Furthermore, the definition of an admissible
element also implies that  $I_1\supset I_2\supset\ldots $  and $I_k$ is the
set of roots of  $\ce^k$.
\\[.6ex]
If $w\in \widehat W$ is admissible, then $I_w$ (resp. $\ce_w$) stands for the
corresponding combinatorial (resp. geometric) Abelian ideal. That is,
\[
   I_w=\{\gamma\in\Delta^+ \mid \delta-\gamma \in\widehat N(w)\}
\ \textrm{ and }\  \ce_w=\oplus_{\ap\in I_w}\g_\ap \ .
\]
Conversely, given $I\in \AD$, we write $w\langle I\rangle$ for the
respective admissible element.  Notice that
\[
  \dim\ce_w =\# (I_w) \textrm{ and } l(w)=\sum_{k\ge 1}\dim(\ce_w)^k \ .
\]
Throughout the paper, $I$ or $I_w$ stands for a combinatorial
{\sf ad}-nilpotent ideal. Whenever we wish to stress that $\AD$ depends on
$\be$ and/or $\g$, we write $\AD(\be)$ or $\AD(\g)$ or even
$\AD(\be,\g)$.

\section{The generators of {\sf ad}-nilpotent ideals}  \label{geners}
\setcounter{equation}{0}

\noindent
Let $I$ be an {\sf ad}-nilpotent ideal. We say that $\gamma\in I$
is a {\it generator} of $I$, if $\gamma-\alpha\not\in I$ for all
$\ap\in\Delta^+$.
Obviously, this is equivalent to the fact that $I\setminus\{\gamma\}$
is still an {\sf ad}-nilpotent ideal. Conversely, if $\varkappa$ is a maximal
element
of $\Delta^+\setminus I$ (i.e., $(\varkappa+\Delta^+)\cap \Delta\subset I$),
%and $(\varkappa + I)\cap \Delta =\varnothing$,
then $I\cup\{\varkappa\}$ is an {\sf ad}-nilpotent ideal.
These two procedures show that the following is true.

\begin{s}{Proposition}  \label{ranked}
Suppose $I\subset J$ are two {\sf ad}-nilpotent ideals. Then there is
a chain of ideals $I=I_0\subset I_1\subset\ldots\subset I_m=J$ such that
$\#(I_{i+1})=\#(I_i)+1$. In other words, $\AD$ is a ranked poset, with
cardinality (dimension) of an ideal as the rank function.
\end{s}%
In the geometric setting, the set of generators
has the following description. For an ideal $\ce=\oplus_{\gamma\in I}\g_\gamma
\subset \be$, there is a unique
$\te$-stable space $\tilde\ce\subset \ce$ such that
$\ce=[\be,\ce]\oplus\tilde\ce$. Then $\gamma$ is a generator of
$I$ if and only if it is a root of $\tilde\ce$.
Write $\Gamma(I)$ for the set of generators of $I$.
It is clear that a subset $\Gamma=\{\gamma_1,\dots,\gamma_l\}\subset
\Delta^+$ is the set of
generators for some ideal if and only if $\gamma_i-\gamma_j\not\in Q^+$
for all $i,j$. This means that $\Gamma\subset\Delta^+$ is the set of
generators for some
{\sf ad}-nilpotent ideal if and only if it is an {\it antichain\/} of $(\Delta^+,\cyeq)$.
This is a manifestation of a general fact that, for any poset $P$,
there is a canonical bijection between the antichains and the dual
order ideals of $P$, see e.g. \cite[3.1]{rstan0}.
\\[.6ex]
In what follows, we also write $I(\Gamma)$
(resp. $\ce(\Gamma)$) for the combinatorial (resp. geometric)
{\sf ad}-nilpotent ideal with the set of generators $\Gamma$.
For instance, the unique maximal element of $\AD$ has the following presentation:
\[
  \begin{array}{rl} \textrm{Geometric: } &   \quad
  \ce(\Pi)=[\be,\be] =\ut \ ; \\
  \textrm{Combinatorial: } & \quad I(\Pi)=\Delta^+ \ .
  \end{array}
\]
It is helpful to have a description of the generators of $I$ in terms of the
respective admissible element.
As usual, we write $\gamma>0$ (resp. $\gamma<0$), if
$\gamma\in\widehat\Delta^+$ (resp. $\gamma\in -\widehat \Delta^+$).

\begin{s}{Theorem}  \label{gener}
Suppose $\gamma\in I_w$. Then $\gamma$ is a generator of $I_w$ if and only if
$w(\delta-\gamma)\in-\widehat \Pi$.
\end{s}\begin{proof*}
``$\Leftarrow$".
Suppose $\gamma$ is not a generator of $I_w$, i.e., $\gamma=\bar\gamma+\nu$,
where $\bar\gamma\in I_w$ and $\nu\in\Delta^+$.
Then $w(\delta-\gamma)=w(\delta-\bar\gamma)-w(\nu)$ is the sum of two
{\sl negative\/} roots.

``$\Rightarrow$".
Set $w(\delta-\gamma)=-\mu<0$. If $\mu$ is not simple, then
$\mu=\mu_1+\mu_2$, where both summands are positive.
%$w(\delta-\gamma)=-\kappa_1-\kappa_2$, where $\kappa_i\in\widehat\Delta^+$.
We have
$w^{-1}(\mu_1)+w^{-1}(\mu_2)=-(\delta-\gamma)< 0$.
Assume for definiteness that $w^{-1}(\mu_2)< 0$.
Since $w^{-1}(-\mu_2)>0$ and $w(w^{-1}(-\mu_2))< 0$, we
have $w^{-1}(-\mu_2)\in \widehat N(w)$, i.e.,
$w^{-1}(\mu_2)=-k\delta+\nu$, where $k\ge 1$ and $\nu\in I_w\subset\Delta^+$.

(a) $k=1$.\\
%for some $\gamma_2\in I_w\subset\Delta^+$.
It follows that
$w^{-1}(\mu_1)=\gamma-\nu\in\Delta$ and $w(\nu-\gamma)=-\mu_1<0$. Since
$w$ is admissible, $\nu-\gamma$ must be negative, i.e.,
$\gamma-\nu\in\Delta^+$. This means that $\gamma$ is not a generator of $I_w$.

(b) $k\ge 2$. \\
Let us show that there is another decomposition of $\mu$ as a sum of
two positive roots such that one has
$k=1$ for one of the summands.
We argue by induction on $k$.
\\[.6ex]
Since $w(k\delta-\nu)<0$, we have $\nu\in (I_w)_k$.
Therefore there is a decomposition
$k\delta-\nu=k'\delta-\nu' + k''\delta-\nu''$, where $k',k''>0$ and
$\nu',\nu''\in I_w$. Hence $\mu_2=\mu_2'+\mu_2''$, where
$w^{-1}(\mu'_2)=\nu'-k'\delta$ and $w^{-1}(\mu''_2)=\nu''-k''\delta$.
The following lemma shows that, in this situation,
$\mu'_2+\mu_1\in\Delta^+$ or  $\mu''_2+\mu_1\in\Delta^+$.
If the latter holds,
then $\mu=\mu'_2 + (\mu''_2+\mu_1)$ is a decomposition such that
$w^{-1}(\mu'_2)=\nu'-k'\delta$, and $k'< k$. This completes the induction step.
\end{proof*}%

\begin{s}{Lemma}  \label{3-kornya}
Suppose $\mu_1,\mu_2,\mu_3\in \widehat\Delta^+$ and
$\mu:=\mu_1+\mu_2+\mu_3\in \widehat\Delta^+$. Then
$\mu_1+\mu_2\in \widehat\Delta^+$ or $\mu_1+\mu_3 \in \widehat\Delta^+$.
\end{s}\begin{proof*}
If $(\mu_2+\mu_3,\mu_1)<0$, then $(\mu_2,\mu_1)<0$ or
$(\mu_3,\mu_1)<0$, and we are done.
If $(\mu_2+\mu_3,\mu_1)\ge 0$, then $(\mu_2+\mu_3,\mu)>0$. Hence
$\mu-\mu_2\in \widehat\Delta$ or $\mu-\mu_3\in \widehat\Delta$.
\end{proof*}%
\begin{s}{Corollary}  \label{number}
The number of generators of\/ $I_w$ is equal to the number of
roots $\ap\in\widehat\Pi$ such that $w^{-1}(\ap)<0$.
\end{s}\begin{proof}
By the definition of an admissible element, if
$w^{-1}(\ap)<0$, then $w^{-1}(\ap)=\gamma-\delta$ for some
$\gamma\in\Delta^+$. Hence $w(\delta-\gamma)\in -\widehat\Pi$
and $\gamma$ is a generator of $I_w$. The rest is clear.
\end{proof}%
Thus, the set of generators of $I_w$ corresponds to a certain subset of
$\widehat \Pi$. More precisely, if $w(\gamma-\delta)=\ap\in\widehat \Pi$
($\gamma\in\Delta^+$),
then we say that $\gamma$ is the generator of $I_w$
{\it corresponding to\/} $\ap$.
\\[.6ex]
Recall that the {\it class of nilpotence\/} of $I\in\AD$, denoted
$\cl(I)$, is the maximal $k$
such that $I_k\ne\varnothing$.
Making use of the admissible element $w$ defining the
{\sf ad}-nilpotent ideal $I_w$, one can readily determine $\cl(I_w)$.

\begin{s}{Proposition}   \label{ind-nil}
$\cl(I_w)=k$ if and only if
$w(\ap_0)+k\delta \in \Delta^+\cup (\delta-\Delta^+)$.
\end{s}\begin{proof}
Since each $(I_w)_m$ is an {\sf ad}-nilpotent ideal, we have
$(I_w)_m\ne \varnothing$ if and only if $\theta\in (I_w)_m$.
Therefore, the very definition of the admissible element
corresponding to an {\sf ad}-nilpotent ideal (see \re{ideals})
implies that
$\cl(I_w)=k$ if and only if
$w(k\delta-\theta)<0$ and $w((k+1)\delta-\theta)>0$.
In other words,
$w(\ap_0)+(k-1)\delta <0$ and $w(\ap_0)+k\delta>0$. Hence the conclusion.
\end{proof}%
{\sf Remark.} If $I$ is a non-trivial Abelian ideal, then $\cl(I)=1$ and the
Proposition asserts that $w(\ap_0)+\delta\in \Delta^+\cup (\delta-\Delta^+)$.
However, Proposition~2.4 in \cite{long} says that
only the first possibility actually realizes,
i.e., $w(\ap_0)+\delta\in \Delta^+$. But, it can be shown that in
case $k=\cl(I)>1$ we do have
both possibilities for $w(\ap_0)+k\delta$.

\begin{rem}{Example}         \label{heisenberg}
Take $w=s_\theta s_0\in \widehat W$, where $s_\theta\in W$ is the reflection with
respect to $\theta$. Then
\[
  s_\theta s_0(\ap)=\left\{ \begin{array}{rc}
  \ap+\delta, & (\ap,\theta)\ne 0 \\ \ap, & (\ap,\theta)=0
  \end{array}\right.  \textrm{ for }\ \ap\in\Pi \ .
\]
We also have
\[
   w^{-1} : \left\{\begin{array}{ccll}
  \ap_0 & \mapsto &\ap_0+2\delta, \\
  \ap_i & \mapsto &\ap_i & \textrm{if } (\ap_i,\theta)=0,\ i\ne 0,\\
  \ap_i & \mapsto &\ap_i-\delta  & \textrm{if } (\ap_i,\theta)\ne 0,\
i\ne 0 \ .
   \end{array}\right.
\]
Hence $w$ is admissible. The corresponding combinatorial
{\sf ad}-nilpotent ideal is $\ch=\{\gamma\in\Delta^+\mid (\gamma,\theta)>0\}$
and the set of generators is $\Gamma(\ch)=\ch\cap\Pi$. The
(geometric) ideal $\ce=\oplus_{\gamma\in\ch}\g_\gamma$ is the standard
Heisenberg subalgebra of $\g$. Obviously, $\cl(\ch)=2$, and we have
$s_\theta s_0(\ap_0)+2\delta=\delta-\theta$.
\end{rem}%
The work of Cellini and Papi \cite{CP2} establishes a bijection between
the {\sf ad}-nilpotent ideals of $\be$ and the points of certain simplex
in $V$ lying in $Q^\vee$, the coroot lattice. This
was used for giving a uniform proof of the formula
for the number of {\sf ad}-nilpotent ideals.
Below, we describe that bijection in a form adapted to our notation,
and show that this can also be used for determining the number
of generators of an ideal.
\\[.6ex]
As is well known, $\widehat W$ is isomorphic to a semi-direct product
of $W$ and $Q^\vee$.
Given $w\in\widehat W$, there is a unique decomposition
\begin{equation}  \label{affine}
w=v_w{\cdot}t_{r_w}\ ,
\end{equation}
where $v_w\in W$ and  $t_{r_w}$ is the translation
corresponding to $r_w\in Q^\vee$.
The word ``translation" means the following.   The group $\widehat W$ has two natural
actions:

(a) the linear action on $\widehat V=V\oplus{\Bbb Q}\delta\oplus{\Bbb Q}\lb$;

(b) the affine-linear action on $V$.
\\
%%We use `$\circ$' to denote the second action.
For $r\in Q^\vee$, the linear action of $t_r\in \widehat W$ on
$V\oplus {\Bbb Q}\delta$ is given by $t_r(x)=x-(x,r)\delta$ (we do not need
the formulas
for the whole of $\widehat V$), while the affine-linear action on $V$ is given
by $t_r\circ y=y+r$. So that $t_r$ is a true translation for this action on $V$.
For instance, the formulae of Example~\ref{heisenberg} show that
$s_\theta s_0=t_{-\theta^\vee}$.
\\
There is a simple procedure for obtaining the affine-linear action on $V$ from the linear
action on $\widehat V$, which is explained in \cite{CP2}, but we do not need this.
\\[.6ex]
%(i.e., $t_r(x)=x+t_r$ for $x\in V$ and $r\in Q^\vee$).
Using the decomposition \re{affine}, one can define the mapping
$\widehat W\to Q^\vee$ by $w\mapsto v_w(r_w)=:d_w$.
One of the main results of \cite{CP2} is that
the set $\{d_w\}$, where $w$ ranges over all admissible elements of
$\widehat W$,  provides a nice parametrization of {\sf ad}-nilpotent ideals.
Namely, set
\[
  \widetilde D=\{ \tau\in V \mid (\tau,\ap)\ge -1\ \forall\ap\in\Pi\ \ \&
  \ \ (\tau,\theta)\le 2\} .
\]
It is a simplex in $V$.
The following is Proposition~3 in \cite{CP2}.

\begin{s}{Theorem}{\ququ (Cellini--Papi)} \ \label{cp2}
The mapping $\AD \to Q^\vee$, defined by
$I\mapsto w\langle I\rangle\mapsto d_{w\langle I\rangle}=:d_I$,
sets up a bijection between $\AD$ and $\widetilde D\cap Q^\vee$.
\end{s}%
{\sf Remark.} Our $\widetilde D\cap Q^\vee$ is $-D$ in the notation
of \cite{CP2}.
\\
Now, we provide a link between the number of generators of
$I$ and the position
of $d_I$ inside of $\widetilde D$.
\begin{s}{Theorem}     \label{main2}
The number of generators of $I$ equals the codimension (in $V$) of the minimal
face of $\widetilde D$ containing $d_I$.
\end{s}\begin{proof}    We have $w=w\langle I\rangle$ is an admissible element
of $\widehat W$.   Let us realise how the vector $d_I=v_w(r_w)$
can be determined by the linear action of $w$.
If $w=v_w{\cdot}t_{r_w}$, then
$w^{-1}=v_w^{-1}{\cdot}t_{-v_w(r_w)}$.    In the following computations, we
repeatedly use the facts that  $\delta$ is isotropic and $w(\delta)=\delta$
for all $w\in
\widehat W$.  If $x\in V\oplus {\Bbb Q}\delta$, then
\[
  w^{-1}(x)=v_w^{-1}(t_{-v_w(r_w)}(x))=
  v_w^{-1}(x+(x,v_w(r_w))\delta)=
  v_w^{-1}(x)+(x, v_w(r_w))\delta=v_w^{-1}(x)+(x, d_I)\delta \ .
\]
In particular, we have
\[
  w^{-1}(\ap_i)=v_w^{-1}(\ap_i)+(\ap_i, d_I)\delta, \quad i\ge 1 ,
\]
and
\[
  w^{-1}(\ap_0)=v_w^{-1}(\ap_0)+(\ap_0, d_I)\delta=-v_w^{-1}(\theta)+
  (1-(\theta,d_I))\delta \ .
\]
Note that $v_w^{-1}(\ap_i)$ and $-v_w^{-1}(\theta)$ are in $\Delta$. Therefore,
by the very definition of an admissible element, we have $(\ap_i, d_I)\ge -1$ \
($i\ge 1$) and
$1-(\theta,d_I)\ge -1$, i.e., $(\theta,d_i)\le 2$.
(In particular, we have recovered the fact that $d_I\in \tilde D$.)
Set $k_i=(\ap_i,d_I)$ and $k_0=1-(\theta,d_I)$.
By Theorem~\ref{gener},
we have $k_i=-1$ if and only if $v_w^{-1}(\ap_i)$ is a generator of $I$; that is,
$I$ has a generator corresponding to $\ap_i$. Similarly, $k_0=-1$ if and only if
$I$ has a generator corresponding to $\ap_0$.
It remains to observe that $k_i=-1$ ($i=0,1,\ldots,p$) are precisely the equations
of facets of $\tilde D$.
\end{proof}%
It follows that an {\sf  ad}-nilpotent ideal has at most $n$
generators%%%%. This fact can be proved directly.
, and the ideals having exactly $n$ generators correspond to
the {\it integral\/} (i.e., lying in $Q^\vee$) vertices of $\tilde D$.
Next, we give an elementary proof for the first observation and
show that $\tilde D$ always has a unique integral vertex.

\begin{s}{Proposition}  \label{vershiny}
Let \ $\Gamma\subset \Delta^+$ be an antichain. Then
\begin{itemize}
\item[\sf (i)] \ The elements of \ $\Gamma$ are linearly
independent and hence $\#(\Gamma)\le\rk\g$;
\item[\sf (ii)] \ If \ $\#(\Gamma)=\rk\g$, then \ $\Gamma=\Pi$.
\end{itemize}
\end{s}\begin{proof*}
(i) Suppose $\Gamma=\{\gamma_1,\dots,\gamma_t\}$.
Since $\gamma_i-\gamma_j\not\in\Delta$, the angle between any pair of
elements of $\Gamma$ is non-acute. Because all $\gamma_i$'s lie in open
half-space of $V$, they are linearly independent.
\\[.5ex]
(ii) Suppose $\Gamma=\{\gamma_1,\dots,\gamma_p\}$, and let $w\in
\widehat W$ be the corresponding admissible element.

We argue by induction on $p$. The case $p=1$ being obvious, we assume that
$p\ge 2$.
If \ $\Gamma\cap\Pi\ne\varnothing$, say $\gamma_1\in\Pi$, then
$\{\gamma_2,\dots,\gamma_p\}$ is an antichain in a root system whose rank is $p-1$.
Hence $\Gamma=\Pi$ by the induction assumption. So, we have only to prove that
the case $\Gamma\cap\Pi=\varnothing$ is impossible. Assume not, i.e.,
$\hot(\gamma_i)\ge 2$ for all $i$. By Theorem~\ref{gener},
\begin{equation}  \label{vyrazheniya}
w(\gamma_i)-\delta
=\ap_{l_i}\in\widehat\Pi .
\end{equation}
%%is a simple root.
Since $p\ge 2$, we may choose $i$ such that $\ap_{l_i}$
lies in $\Pi$. Without loss of generality, we may assume that
$i=1$.
Choose also roots $\mu,\bar\mu\in\Delta^+$
such that $\gamma_1=\mu+\bar\mu$. Obviously, then $\mu,\bar\mu\not\in I=I(\Gamma)$.
By part (i), $\Gamma$ is a basis for $V$.
Hence,
\[
  \mu=\sum_{j\in J}d_j\gamma_j- \sum_{k\in K}c_k\gamma_k \ ,
\]
where $J,K$ are disjoint subsets of $\{1,2,\dots,p\}$ and $c_k,d_j>0$.
Therefore,
\[
  \bar\mu=-\sum_{j\in J}d_j\gamma_j+ \sum_{k\in K}c_k\gamma_k +\gamma_1 \ .
\]
Given $\nu\in \widehat \Delta$, we say that the {\it level\/}
of $\nu$, denoted $\lev(\nu)$, is $m\in\Bbb Z$,
if $\nu-m\delta\in \Delta$.
Consider the roots $w(\mu), w(\bar\mu)\in \widehat \Delta$.
Since $w(\delta-\mu)>0$ and $w(\mu)>0$, we have $\lev(w(\mu))$ is either 1 or
0, and likewise for $\bar\mu$. As  $w(\mu+\bar\mu)=\delta+\gamma_1$ has level 1,
we may assume without loss that $\lev(w(\mu))=0$ and
$\lev(w(\bar\mu))=1$.
Using Eq.~\ref{vyrazheniya} for the $w(\gamma_i)$'s, we obtain
\[
 w(\mu)=(\sum_j d_j-\sum_k c_k)\delta +
\sum_{j\in J}d_j\ap_{l_j}- \sum_{k\in K}c_k\ap_{l_k} \ ,
\]
and
\[
 w(\bar\mu)=(1-\sum_j d_j+\sum_k c_k)\delta -
\sum_{j\in J}d_j\ap_{l_j}+ \sum_{k\in K}c_k\ap_{l_k}+ \ap_{l_1} \ .
\]
If one of the roots $\ap_{l_i}$, $i\in J\cup K$, is equal to $\ap_0=\delta-\theta$,
then the equality $\lev(w(\bar\mu))-\lev(w(\mu))=1$ cannot be satisfied.
Hence all these roots lie in $\Pi$ and hence
$\sum_j d_j-\sum_k c_k=\lev(w(\mu))=0$.
But the equality $w(\mu)=\sum_{j\in J}d_j\ap_{l_j}- \sum_{k\in K}c_k\ap_{l_k}\in
\Delta$
contradicts the fact that $w(\mu)$ is positive.
\end{proof*}%

\begin{s}{Corollary}   \label{odna}
The simplex $\tilde D$ has a unique integral vertex, corresponding to the
unique maximal {\sf ad}-nilpotent ideal.
\end{s}%
The vertices of $\tilde D$ can explicitly be described, see \cite{CP2}.
Indeed, let $\{\pi_i\}$ be the basis for $V$ dual to
$\{\ap_i\}$, $1\le i\le p$ and $h$ the Coxeter number for $\Delta$.
If $\theta=\sum_{i=1}^p m_i\ap_i$ and
$\rho^\vee$ is the half-sum of all positive coroots, then the vertices of
$\tilde D$ are $-\rho^\vee$ and $-\rho^\vee+\frac{h+1}{m_i}\pi_i$, $1\le i\le p$.
However, it is not immediately clear from this that exactly one vertex lies
in $Q^\vee$.

\section{A combinatorial statistic on $\AD(\g)$, Catalan arrangements, and clusters}
\label{q-sign}
\setcounter{equation}{0}

By \cite{CP2}, the cardinality of $\AD(\g)$ is equal to
$\displaystyle \prod_{i=1}^p \frac{h+e_i+1}{e_i+1}$, where the $e_i$'s
are the exponents and $h$ is the Coxeter number of $\g$.
In this section, we consider the {\it simple root statistic\/} on
$\AD(\g)$. It is given by
\[
   \text{sim}(I)=\#(I\cap \Pi),\qquad I\in \AD(\g) \ .
\]
Accordingly, we set
\[
\AD(\g)_i=\{ I\in\AD(\g) \mid \#(I\cap \Pi)=i\},\quad i=0,1,\dots,p \ .
\]
Because we know the number of {\sf ad}-nilpotent ideals for all
simple $\g$, the number $\AD(\g)_0$
can be counted via the inclusion-exclusion principle. Indeed, the ideals
containing $\ap_i\in\Pi$ can be identified with the ideals of the semisimple
subalgebra of $\g$ whose simple roots are $\Pi\setminus\{\ap_i\}$.
%Set
%\[
%    \AD(\g)_0=\{ I\in\AD(\g) \mid  I\cap\Pi= \varnothing \} \ ,
%\]
Write $\g(J)$ for the semisimple subalgebra of $\g$ whose set of simple
roots is $J\subset \Pi$. Then
\begin{equation}  \label{sim-zero}
   \# \AD(\g)_0=\sum_{J\subset \Pi} (-1)^{p-\# J}\# \AD(\g(J)) \ .
\end{equation}
In turn, the numbers $\AD(\g)_i$ $(i>0)$ are easily computed, once one knows
$\AD(\g)_0$. For instance, the number of ideals containing exactly one
simple root, say $\ap_i$, is equal to the number of
all ideals in $\g(\Pi\setminus \{\ap_i\})$ that do not contain
simple roots. Hence
\[
   \#\AD(\g)_1= \sum_{\# J=p-1} \# \AD(\g(J))_0 \ .
\]
Similarly, one obtains the general formula:
\begin{equation}  \label{sim}
   \#\AD(\g)_i= \sum_{\# J=p-i} \# \AD(\g(J))_0 \ .
\end{equation}
Of course, applying Equations~\re{sim-zero} and \re{sim}, one should use
the relation that if $\h=\h_1\oplus \h_2$, then $\#\AD(\h)=\#\AD(\h_1){\cdot}
\#\AD(\h_2)$, and likewise for $\AD(\h)_0$.
\\[.6ex]
The distribution of the simple root statistic over
$\AD(\g)$ yields the polynomial
\[
   \cs_\g(q)=\sum_{i=0}^p \#(\AD(\g)_i) q^i \ ,
\]
which is not hard to compute. For instance, the following table contains the
relevant data for exceptional Lie algebras.

\begin{center}
\begin{tabular}{c|c|rrrrrrrrr}
\multicolumn{2}{c|}{$i$} & 0 & 1 & 2 & 3 & 4 & 5 & 6 & 7 & 8 \\ \hline
   & $\GR{F}{4}$ &   66 &   24 &   10 &   4 &   1  & \\
$\#\AD(\g)_i$
  & $\GR{E}{6}$ &   418 &  228 &  110 &  50 &  20 &   6 &  1 & \\
  & $\GR{E}{7}$ &  2431 & 1001 &  429 & 187 &  77 &  27 &  7 & 1  \\
  & $\GR{E}{8}$ & 17342 & 4784 & 1771 & 728 & 299 & 112 & 35 & 8 & 1 \\ \hline
\end{tabular}
\end{center}
It immediately follows from Eq.~\re{sim} that
$\#\AD(\g)_p=1$ and $\#\AD(\g)_{p-1}=p$.
A bit longer analysis yields
\begin{s}{Proposition}  \label{p-2}

If $\g$ simply-laced, then $\#\AD(\g)_{p-2}=(p-1)(p+2)/2$;

If $\g\in \{B,C,F\}$, then  $\#\AD(\g)_{p-2}=p(p+1)/2$.
\end{s}\begin{proof}
There are $p(p-1)$ subalgebras of the form $\g(J)$ with $\# J=2$.
Of these subalgebras, we have \\
$\bullet$\quad
$p-1$ subalgebra of type $\GR{A}{2}$ and $(p-1)(p-2)/2$ subalgebras of
type $\GR{A}{1}\times \GR{A}{1}$, if $\g$ is simple laced;\\
$\bullet$\quad
one subalgebra of type $\GR{C}{2}$, $p-2$ subalgebras of type $\GR{A}{2}$ and $(p-1)(p-2)/2$ subalgebras of
type $\GR{A}{1}\times \GR{A}{1}$, if $\g$ is doubly laced.
\end{proof}%
Our results and conjectures for the classical series are as follows.

\begin{s}{Theorem}   \label{sim-sl}
$\#\AD(\GR{A}{p})_i=\displaystyle
\frac{i+1}{p+1}\genfrac{(}{)}{0pt}{}{2p-i}{p}$, \quad $i=0,1,,\dots,p$.
\end{s}%
We defer the proof to Section~\ref{sl}.
Arguing by induction on $p$ and using Eq.~\re{sim}, one obtains
$\cs_\g(q)=\cs_{\g^\vee}(q)$, where $\g^\vee$ is the Langlands dual Lie
algebra. The only practical output of this equality is that the simple
root statistic has the same distribution for $\GR{B}{p}$ and $\GR{C}{p}$.
However, we have only conjectural values for $\GR{C}{p}$ and $\GR{D}{p}$,
which are verified for $p\le 8$.
%However, I believe that, at list
%for $\GR{C}{p}$, the result is known to combinatorialists. That is to say,
%it is buried somewhere in the literature.
\begin{s}{Conjecture}     \label{sim-sp}

$\#\AD(\GR{C}{p})_i=\displaystyle
\genfrac{(}{)}{0pt}{}{2p-1-i}{p-1}$, \quad $i=0,1,,\dots,p$.

$\#\AD(\GR{D}{p})_i=\displaystyle
\genfrac{(}{)}{0pt}{}{2p-2-i}{p-2}+ \genfrac{(}{)}{0pt}{}{2p-3-i}{p-2}$,
 \quad $i=1,2,,\dots,p$.
\end{s}%
Notice that the conjecture does not give an expression for
$\#\AD(\GR{D}{p})_0$. As we will see below, the right value for
$\#\AD(\GR{D}{p})_0$ is
$\genfrac{(}{)}{0pt}{}{2p-2}{p-2}+ \genfrac{(}{)}{0pt}{}{2p-3}{p-3}$.

\noindent
Using Eq.~\re{sim-zero}, it is easy to compute $\#\AD(\g)_0$ for any
simple Lie algebra. However, obtaining a closed expression in the classical
case requires some work.
%In this way, we obtain certain numbers for the exceptional Lie algebras.
%But the classical series require more work to obtain a closed expression.
In order to obtain a more conceptual explanation and the closed formula
valid for all $\g$, we use the theory of arrangements.
\begin{rem}{Remark}    \label{sovpad}
Having written up Propositions~\ref{regions} and \ref{no-simple} below,
I found that exactly the same results are obtained in the recent preprint
of C.~Athanasiadis \cite{ath02}.
In this preprint, he gave a conceptual proof of the formula \re{char-cat}
for the characteristic polynomial of the Catalan arrangement.
In fact, Eq.~\re{char-cat} was known for all simple Lie algebras via
case-by-case verification, and this was used in my original argument.
%As an application,
%he derived an explicit formula for $\#\AD(\g)_0$. In my
\end{rem}%
%\\[.5ex]
Recall a bijection between the {\sf ad}-nilpotent ideals
and the regions of the Catalan arrangement that are contained in the
fundamental Weyl chamber. This bijection is due to
Shi \cite[Theorem\,1.4]{shi}, see also \cite[\S\,4]{CP2}.
The {\it Catalan arrangement\/} $\mathrm{Cat}(\Delta)$
is the set of hyperplanes in $V$ having the equations
\[
  (x,\mu)=1,\qquad (x,\mu)=0, \qquad (x,\mu)=-1 \hspace{7ex} (\mu\in\Delta^+)\ .
\]
The {\it regions\/} of an arrangement are the connected components of
the complement in $V$ of the union of all its hyperplanes.
Clearly, $\ccc$ is a union of regions of $\mathrm{Cat}(\Delta)$.
Any region lying in $\ccc$ is said to be
{\it dominant\/}.
The bijection takes an ideal $I\subset \Delta^+$ to the region
\[
 R_I =\{ x\in \ccc \mid  (x,\gamma)>1, \text{ if \ }\gamma\in I \quad \&\quad
           (x,\gamma)<1,  \text{ if \ }\gamma\not\in I
           \}\ .
\]
Obviously, the dominant regions of $\mathrm{Cat}(\Delta)$ are the same as
those for the {\it Shi arrangement\/} $\mathrm{Shi}(\Delta)$.
Here $\mathrm{Shi}(\Delta)$ is the set of hyperplanes in $V$ having the equations
\[
  (x,\mu)=1,\qquad (x,\mu)=0  \hspace{9ex} (\mu\in\Delta^+)\ .
\]
It will be more convenient for us to deal with the arrangement $\mathrm{Cat}(\Delta)$,
since it is $W$-invariant.
A region (of an arrangement) is called {\it bounded\/}, if it is contained in
a sphere about the origin.

\begin{s}{Proposition}   \label{regions}
$I\in\AD(\g)_0$
%%contains no simple roots
if and only if the region $R_I$ is bounded.
\end{s}\begin{proof}
1. Suppose $I\cap\Pi=\varnothing$. Then the definition of $R_I$ shows that
ii is contained in the bounded domain in $\ccc$ given by the inequalities
$(\ap,x)<1$, $\ap\in\Pi$.

2. Suppose $\ap_i\in I\cap \Pi$. Then $I$ also contains all positive roots
whose coefficient of $\ap_i$ in the expression through the simple roots
is positive. Hence for all roots $\gamma$ such that $(\gamma,\vp_i)>0$
we have the constraints
$(x,\gamma)>1$. This means that if $x\in R_I$, then all constraints
are satisfied for $x+a\vp_i$  with any $a\in {\Bbb R}_{\ge 0}$.
Thus, $R_I$ is unbounded.
\end{proof}%
The number of regions and bounded regions of any hyperplane arrangement
can be counted through the use of a striking result of T.\,Zaslavsky.
Let $\chi(\ca, t)$ denote the characteristic polynomial of a
hyperplane arrangement $\ca$ in $V$ (see e.g. \cite{ath96},\,\cite{sagan} for
precise definitions).
%%The following is proved in \cite{zasl}
\begin{s}{Theorem} \label{zaslavsky} {\ququ (Zaslavsky \cite[Sect.\,2]{zasl})}    \\
1. The number of regions into which $\ca$ dissects $V$ equals
$r(\ca)=(-1)^p\chi(\ca,-1)$.
\\[.6ex]
2. The number of bounded regions into which $\ca$ dissects $V$ equals
$b(\ca)=|\chi(\ca,1)|$.
\end{s}%
Recently, Athanasiadis \cite{ath02} found a rather simple case-free proof of
the following formula for the characteristic polynomial of
the Catalan arrangement:
%The characteristic polynomials of Catalan arrangements for
%the classical series are computed by  C.\,Athanasiadis \cite{ath96}.
\begin{equation} \label{char-cat}
  \chi(\mathrm{Cat}(\Delta),t)=\prod_{i=1}^p (t-h-e_i) \ .
\end{equation}
(For the classical series, it was computed earlier in \cite{ath96}).
%%by  C.\,Athanasiadis \cite{ath96}.
%The exceptional root systems can be treated in a case-by-case fashion.
Now, combining the preceding results, we arrive at our goal.

\begin{s}{Proposition}  \label{no-simple}
$\#\AD(\g)_0=\displaystyle \prod_{i=1}^p \frac{h+e_i-1}{e_i+1}$ .
\end{s}\begin{proof}
Since the arrangement $\mathrm{Cat}(\Delta)$ is $W$-invariant, the number
of its bounded regions lying in $\ccc$ is equal to
$\frac{1}{\# W} |\chi(\mathrm{Cat}(\Delta),1)|$. It remains to observe that
$\# W=\prod_i (e_i+1)$.
\end{proof}%
Similarly, using the value $\chi(\mathrm{Cat}(\Delta),-1)$, as
Athanasiadis also did in \cite{ath02}, one obtains the
formula for the number of all {\sf ad}-nilpotent ideals stated at
the beginning of this section. This proof is not so elementary as the proof
of Cellini-Papi \cite{CP2}, for it requires some deep results from the
theory of arrangements.
\\[.6ex]
It is quite interesting that the numbers
$\displaystyle \prod_{i=1}^p \frac{h+e_i+1}{e_i+1}$
and $\displaystyle \prod_{i=1}^p \frac{h+e_i-1}{e_i+1}$
also appear in \cite[Theorem\,1.9 \& Prop.\,3.9]{cluster} as the numbers of all and positive clusters,
respectively. We are not going to discuss the theory of clusters related to
the root systems, referring to that paper for all relevant definitions.
For our current purposes, it suffices to know that clusters are
certain subsets
of $\Delta^+\cup (-\Pi)$. Each cluster is a linearly independent subset
of $V$ having
exactly $p$  elements. A cluster is called {\it positive\/}, if all its
elements are positive roots.
\\[.5ex]
A close relationship between clusters and {\sf ad}-nilpotent ideals is
seen in the following curious fact.
Let $\CL(\g)_i$ denote the set of clusters having exactly $i$
elements from $-\Pi$.

\begin{s}{Theorem}   \label{empir}
One always has the equality
$\# \AD(\g)_i= \# \CL(\g)_i$.
\end{s}\begin{proof} From
Proposition~3.6 in \cite{cluster}, it follows that the numbers
$\CL(\g)_i$, $i=0,1,\dots,p$,  also satisfy the recurrent relations Eq.~\re{sim-zero} and \re{sim}.
\end{proof}%
It is not too brave to suggest that there exists a natural bijection
between clusters and {\sf ad}-nilpotent ideals that takes
$\CL(\g)_i$ to $\AD(\g)_i$ for all $i$.

%
%Some particular cases of \re{empir} are easy to handle.
%For instance, it is clear that $\# \AD(\g)_p= \# \CL(\g)_p=1$
%and $\# \AD(\g)_{p-1}= \# \CL(\g)_{p-1}=p$.
%\noindent
%One more connection (i.e., numerical coincidence)
%of {\sf ad}-nilpotent ideals with clusters will be distinguished in
%Section~\ref{vmeste}.

\section{On {\sf ad}-nilpotent ideals for $\g={\frak sl}_n$}
\label{sl}
\setcounter{equation}{0}

\begin{figure}[htb]
\setlength{\unitlength}{0.018in}
\begin{center}
\begin{picture}(125,125)(1,0)
\multiput(0,0)(120,0){2}{\line(0,1){120}}
\multiput(0,0)(0,120){2}{\line(1,0){120}}
\put(50,110){\line(1,0){70}}
\put(60,100){\line(1,0){60}}
\put(70,90){\line(1,0){50}}  % \put(95,105){\vector(0,-1){9}}
\put(90,80){\line(1,0){30}}  % \put(105,95){\vector(-1,0){9}}
\multiput(100,40)(0,10){4}{\line(1,0){20}}

\qbezier[40](-5,115),(30,115),(55,115)  % horizontal
\qbezier[77](-5,45),(50,45),(105,45)     % horizontal
\qbezier[50](105,125),(105,85),(105,45)
\qbezier[7](55,125),(55,120),(55,115)

%\put(50,70){\dashbox{1}(0,40){}}
%\put(50,70){\dashbox{1}(70,0){}}
%\put(40,80){\dashbox{1}(0,40){}}
%\put(40,80){\dashbox{1}(50,0){}}

\put(50,110){\line(0,1){10}}
\put(60,100){\line(0,1){20}}
\put(70,90){\line(0,1){30}}
\put(80,90){\line(0,1){30}}
\put(90,80){\line(0,1){40}}
\put(100,40){\line(0,1){80}}
\put(110,40){\line(0,1){80}}

\put(-11,115){$i_1$}
\put(-11,45){$i_k$}
\put(52,127){$j_1$}
\put(75,127){$\cdots$}
\put(102,127){$j_k$}

\put(120,0){\line(-1,1){120}}   %diagonal

\put(1,119){\line(1,0){48}}             \put(49,109){\line(0,1){10}}
\put(49,109){\line(1,0){10}}           \put(59,99){\line(0,1){10}}
\put(59,99){\line(1,0){10}}             \put(69,89){\line(0,1){10}}
\put(69,89){\line(1,0){20}}             \put(89,79){\line(0,1){10}}
\put(89,79){\line(1,0){10}}             \put(99,39){\line(0,1){40}}
\put(99,39){\line(1,0){20}}             \put(119,1){\line(0,1){38}}

\end{picture}
\end{center}
\caption{An {\sf ad}-nilpotent ideal in $\sln$}   \label{pikcha_A}
\end{figure}
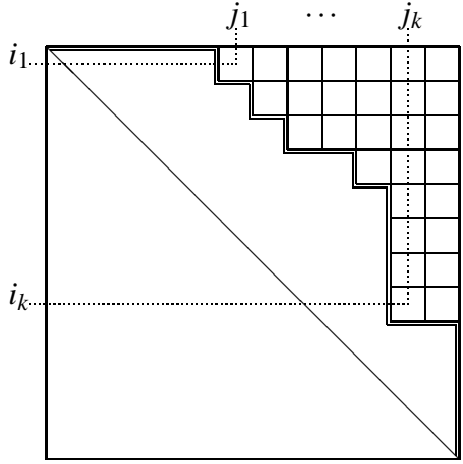

\noindent
At the rest of the paper, we are going to study another combinatorial
statistic on the set of {\sf ad}-nilpotent ideals, which is related to the
theory developed in Section~\ref{geners}.
We first consider the classical series in Sections~\ref{sl} and
\ref{classical}, and then move to the general case in
Section~\ref{vmeste}.
\\[.6ex]
At the rest of this section, $\g=\sln$ and hence $p=n-1$.
We assume that $\be$ (resp. $\te$) is standard, i.e., it is the space of
upper-triangular (resp. diagonal)
matrices. Then the positive roots are identified
with the pairs $(i,j)$, where $1\le i<j\le n$. For instance, $\ap_i=(i,i+1)$ and
$\theta=(1,n)$. An {\sf ad}-nilpotent $\be$-ideal is represented by a
right-justified
Ferrers diagram with at most $n-1$ rows, where the length of $i$-th row
is at most $n-i$. If a box of a Ferrers diagram corresponds to a
positive root $(i,j)$, then we say that this box has the coordinates
$(i.j)$.
The unique northeast corner of the diagram corresponds to $\theta$
and the southwest corners give rise to the
generators of the corresponding ideal, see  Figure~\ref{pikcha_A}.
%We write $\AD_n$ for the set of

\noindent
Such a diagram (ideal) $I$ is completely determined by the coordinates of boxes
%(roots)
that contain the southwest  corners of the diagram, say $(i_1,j_1),\dots,
(i_k,j_k)$. Then
we obviously  have  $\Gamma(I)=\{(i_1,j_1),\dots, (i_k,j_k)\}$ and
\[  %begin{equation}  \label{ogran}
1\le i_1< i_2< \dots  <i_k\le n-1,\quad
  2\le j_1< j_2< \dots  <j_k\le n \ .
\]  %end{equation}
Various enumerative results for {\sf ad}-nilpotent ideals in $\sln$ are obtained
in \cite{AKOP},\,\cite{CP1},\,\cite{OP}.
In particular, the total number of {\sf ad}-nilpotent ideals equals
$C_n=\frac{1}{n+1}\genfrac{(}{)}{0pt}{}{2n}{n}$, the $n$-th Catalan number.
There is a host of combinatorial objects that are counted by Catalan
numbers, see \cite[ch.\,6, Ex.\,6.19]{rstan} and the ``Catalan addendum'' at
\verb|www-math.mit.edu/~rstan/ec|.
We shall use the fact that $C_n$ is equal to

(a) the number of   all sequences
$v=v_1v_2\ldots v_{2n}$ of $n$ 1's and $n$ \ $-1$'s with all partial sums
nonnegative,
or   %%  Ex. 19(r)

(b) the number of lattice paths from $(0,0)$ to $(n,n)$ with steps $(1,0)$ and $(0,1)$,
%never rising above   the line $x=y$.
always staying in the domain $x\le y$, i.e., the number of Dyck paths of
semilength $n$.

In our matrix interpretation, we are forced to assume that the
$x$-axis is vertical and directed downstairs, while the $y$-axis is horizontal.
Therefore $(0,0)$ is the upper-left corner and
$(n,n)$ is the lower-right corner of the matrix.
The Dyck path corresponding to an
{\sf ad}-nilpotent ideal is the double path in Figure~\ref{pikcha_A}.
It has $2n$ steps. The corresponding sequence $v$ is obtained as follows.
%The lattice path consists of $2n$ steps;
We start from
$(0,0)$ and attach  $+1$ to the horizontal step (i.e., $(0,1)$) and
$-1$ to the vertical step (i.e., $(1,0)$).

{\bf Remark.} Coordinates of boxes of Ferrers diagrams and lattice points
considered above are compatible in the sense that the coordinates of
a box are equal to the coordinates of its southeast corner.

{\sl Proof of Theorem~\ref{sim-sl}.}  Once the relationship between the
{\sf ad}-nilpotent ideals in $\sln$ and Dyck paths is established, one
may appeal to huge combinatorial literature on the latter.
It is clear that $I\in\AD(\sln)$ contains a simple root if and only if
the corresponding Dyck path touches the diagonal somewhere except
the points $(0,0)$ and $(n,n)$. In other words, the number of
simple roots in $I$ equals the number of (intermediate) {\it returns\/} \
of the Dyck path. The distribution of this statistic
%%on the Dyck paths
is well-known, see e.g. \cite[6.6]{deu}.  \qus
Let $\AD_n$ denote the set of all {\sf ad}-nilpotent ideals for $\sln$.
{From} now on, we stick to considering the statistic \
$\text{gen}: \AD_n \to {\Bbb N}$,
which assigns to an ideal the number of its generators.
Let
$\AD_{n}^{k}$, $0\le k\le n{-}1$,
be the set of  ideals with $k$ generators, i.e.,
the set of Ferrers diagrams, as above, with exactly $k$ southwest corners.

\begin{s}{Proposition}    \label{narayana}
$\# (\AD_{n}^{k})=\displaystyle \frac{1}{n}\genfrac{(}{)}{0pt}{}{n}{k}
\genfrac{(}{)}{0pt}{}{n}{k+1}$.
\end{s}\begin{proof}
The numbers $N(n,k)=\frac{1}{n}\genfrac{(}{)}{0pt}{}{n}{k}
\genfrac{(}{)}{0pt}{}{n}{k-1}$, $k=1,\dots,n$,
are called the {\it Narayana numbers}, so that we are to show that
$\#(\AD_{n}^{k-1})=N(n,k)$. It is known that the Narayana numbers have the following
combinatorial interpretation, see \cite[ch.\,6, Ex.\,36(a)]{rstan}.
Let $X_{nk}$ be the set of all sequences $v=v_1\dots v_{2n}$ as in (a) above,
%$w=w_1w_2\ldots w_{2n}$ of $n$ 1's and $n$ \ $-1$'s with all partial sums nonnegative,
such that
\[
  k=\# \{j \mid v_j=1, v_{j+1}=-1\} \ .
\]
Then $\#(X_{nk})=N(n,k)$.
For $v\in X_{nk}$, it is easily seen that
\[
  k-1=   \# \{j \mid v_j=-1, v_{j+1}=1\} \ .
\]
The change of sign from 1 to $-1$ (resp. from $-1$ to 1)
in $v$ corresponds to the turn of the type ``horizontal followed by vertical''
(resp. ``vertical followed by horizontal'') step in the respective lattice path.
Geometrically, the steps of second type correspond to the southwest corners
of our Ferrers diagram.
It follows that the sequences $v\in X_{nk}$ are in bijection with the
Ferrers diagrams with $k-1$ southwest corners, and we are done.
%Indeed, the Catalan number $C_n$ counts the number of   Dyck paths from
%$(0,0)$ to $(2n,0)$ never falling below the $x$-axis
%\cite[ch.\,6, Ex.\,19(i)]{rstan}. Then $X_{nk}$
%corresponds  to Dyck path with exactly $k$ peaks (maxima). Clearly,
%there are $k-1$ minima between these maxima.
%It remains to notice that the minima of the Dyck path corresponding to
%an {\sf ad}-nilpotent ideal are in one-to-one correspondence with the generators of
%the ideal.
\end{proof}%
Since $N(n,k)=N(n,n-k+1)$, one may suggest that
there is a bijective interpretation of this equality. This is really the case.

\begin{s}{Theorem}   \label{1-1}
There is a natural bijection between
%the sets
$\AD_{n}^{k}$ and $\AD_{n}^{n-k-1}$.
\end{s}\begin{proof}
Let $(i_1,j_1),\dots,(i_k,j_k)$ be the generators of an {\sf ad}-nilpotent ideal
$I\in \AD_{n}^{k}$. Consider separately the ordered sets
of the first and second coordinates for these
generators, i.e.
put $X(I)=\{i_1,\dots,i_k\}$ and $Y(I)=\{j_i,\dots,j_k\}$.
We wish to construct two other ordered sets
that will form the first and the second  coordinates of the generators
for the dual ideal. To this end, put

$X(I^*)=\{1,\dots,n{-}1\}\setminus \{j_1-1,\dots,j_k-1\}$.

$Y(I^*)=\{2,\dots,n\}\setminus \{i_1+1,\dots,i_k+1\}$.
\\
For $A=\{a_1,\dots,a_m\}$, it is convenient to
introduce notation $A[a]=\{a_1+a,\dots,a_m+a\}$. Then the previous formulas
can be written as
\begin{equation} \label{dual}
\begin{array}{l}
X(I^*)=(\{2,\dots,n\}\setminus Y(I))[-1], \\
Y(I^*)=(\{1,\dots,n-1\}\setminus X(I))[1].
\end{array}
\end{equation}
It is then easily seen that the square of this transformation
is the identity on $\AD_{n}^{k}$. Therefore one has only to prove that the
ordered sets $X(I^*),Y(I^*)$ determine an {\sf ad}-nilpotent ideal.
The latter means that if  $X(I^*)=\{i_1^*,\dots,i_{n-k-1}^*\}$ and
$Y(I^*)=\{j_1^*,\dots,j_{n-k-1}^*\}$, then $i^*_q < j^*_q$ for all
$q$. (Of course, $i^*_1<i^*_2<\dots $ and likewise for $j^*_l$.)

(a) Given $q\in \{1,\ldots,n-k-1\}$, suppose there is $m$ such that
$i_m>m+q-1$. Assume also that $m$ is the minimal number with this
property. Then $i_m\ge m+q$ and $i_{m-1}< m-1+q$. Therefore the
$q$-th element of $\{1,\dots,n-1\}\setminus X(I)$ is $m-1+q$ and
hence $j^*_q=m+q$. Since $j_m>i_m=m+q$, we can find the minimal
number $l$ such that $j_l> l+q$. Then $l\le m$ and the $q$-th element
of $\{2,\dots,n\}\setminus Y(I)$ is $l+q$.
Thus, $i^*_q=l+q-1 < m+q=j^*_q$.

(b) Suppose $i_m\le m+q-1$ for all $m\in\{1,\dots,k\}$, that is,
$i_k\le k+q-1$. Then the $q$-th element of $\{1,\dots,n-1\}\setminus X(I)$ is
$k+q$ and hence $j^*_q=k+q+1$. On the other hand, the inequalities
$i^*_q<i^*_{q+1}< \dots i^*_{n-k-1}\le n-1$ show that
$i^*_q\le (n-1)-((n-k-1)-q)=q+k$.
\\[.7ex]
Thus, $X(I^*)$ and $Y(I^*)$ determine an element of $\AD_{n}^{n-k-1}$,
which we denote by $I^*$.
\end{proof}%
For all $k\in\{0,1,\dots,n{-}1\}$, we have constructed  bijections
\[
     \AD_{n}^{k}\to \AD_{n}^{n-k-1}, \quad I\mapsto I^* \ .
\]
which give rise to an involutory transformation $\ast: \AD_n\to \AD_n $.
Although this transformation is not order-reversing with respect to the
inclusion of ideals, it has interesting properties.
%Let $h=h(\g)$ denote the Coxeter number of $\g$.
The formulation of these properties
is ``universal'', i.e., it makes sense for any (semi)simple Lie algebra:
\begin{s}{Lemma}  \label{universal}
Suppose $A\subset\Pi$ is an arbitrary subset, and  $I=I(A)$.
Then $I^*=I(\Pi\setminus A)$.
\end{s}\begin{proof}
Straightforward. Use Formulae~\re{dual}.
\end{proof}%
To state one more property, we need some notation.
As usual, the height of a root $\gamma\in\Delta^+$ is denoted by $\hot(\gamma)$.
Recall that $h=\hot(\theta)+1$ is the Coxeter number of $\g$.
Set $\Delta^+(k)=\{\gamma \in \Delta^+\mid \hot(\gamma)=k\}$ and
$\Delta^+_k=\{\gamma \in \Delta^+\mid \hot(\gamma)\ge k\}$.
It is clear that $\Delta^+_k$ is a combinatorial {\sf ad}-nilpotent ideal and
$\Gamma(\Delta^+_k)=\Delta^+(k)$.
\\[.6ex]
For $\sln$, we have $\hot(i,j)=j-i$ and the Coxeter number is $n$.

\begin{s}{Lemma}  \label{universal2}
In case of\/ $\sln$, we have
$(\Delta^+_k)^*=\Delta^+_{h{+}1{-}k}=\Delta^+_{n{+}1{-}k}$.
\end{s}\begin{proof}   Set $I=\Delta^+_k$.
In our notation, the roots in $ \Delta^+(k)$ are $(1,k+1),(2,k+2),\dots, (n-k,n)$.
Hence $X(I)=\{1,2,\dots,n-k\}$ and  $Y(I)=\{k+1,k+2,\dots,n\}$.
Therefore $X(I^*)=\{1,2,\dots,k-1\}$ and $Y(I^*)=\{n-k+2,\dots,n\}$.
This means that $I^*$ is generated by the roots $(1,n-k+2),\dots,(k-1,n)$, i.e.,
all roots of height $n-k+1$.
\end{proof}%
{\sf Examples.} In the geometric context,
taking $k=1$, we obtain $\ut^*=\{0\}$.
For $k=2$, we have $[\ut,\ut]^*=\g_\theta$, because $\theta$ is the only root
of height $h-1$.
\\[.6ex]
It is curious that our definition of the dual {\sf ad}-nilpotent ideal
for $\sln$ leads to another occurrence of Catalan numbers. Namely,
let us try to describe and enumerate the self-dual ideals.
For $I\in \AD_{n}^{m}$, the necessary condition of self-duality is
$m=n-m-1$. That is, $n=2m+1$.

\begin{s}{Theorem}  \label{self-dual}
There are no self-dual\/ {\sf ad}-nilpotent ideals for\/ ${\frak sl}_{2m}$.
For\/ ${\frak sl}_{2m+1}$, the number of self-dual\/ {\sf ad}-nilpotent
$\be$-ideals is equal to \
$\displaystyle  \frac{1}{m+1}\genfrac{(}{)}{0pt}{}{2m}{m}$.
\end{s}\begin{proof}
We use the notation introduced in Theorem~\ref{1-1}.
Suppose $I\in \AD_{2m+1}^{m}$ and
$X=X(I)=\{i_1,i_2,\dots,i_m\}$,  $Y=Y(I)=\{j_1,j_2,\dots,j_m\}$.
The condition $I=I^*$ means that
$X=X^*=\bar Y[-1]$ %=(\{2,3,\dots,2m+1\}\setminus Y)[-1]$
and
$Y=Y^*=\bar X[1]$. %=(\{1,2,\dots,2m\}\setminus X)[1]$.
Clearly, all these equalities are equivalent to the following
%That is,
%$\{1,2,\dots,2m\}\setminus Y[-1]=X$.
\[
  \{1,2,\dots,2m\}=\{i_1,i_2,\dots,i_m\} \sqcup   \{j_1-1,j_2-1,\dots,j_m-1\}=
  X \sqcup Y[-1] \ .
\]
Therefore $Y$ is determined by $X$ and vice versa.
However, $X$ cannot be an arbitrary
$m$-element subset of $\{1,2,\dots,2m\}$, since the conditions
$i_k< j_k$, $k=1,\dots,m$, must also be satisfied.
Given $X\subset \{1,2,\dots,2m\}$ with  $\#(X)=m$, define the sequence
$v=v_1v_2\ldots v_{2m}$ by the following rule:
\[
  v_i= \left\{\begin{array}{rc}  1, \textrm{ if \ } i\in X ,\\
             -1, \textrm{ if \ } i\not\in X .
                       \end{array}\right.
\]
Then   $(X, \bar X[1])$ determine an {\sf ad}-nilpotent ideal if and only if
all partial sums of $v$ are nonnegative.
Indeed, $\sum_{i=1}^{2k-1}v_i <0$ if and only if $i_k\ge j_k$.
As was mentioned above, the number
of such sequences is the $m$-th Catalan number.
\end{proof}%
To illustrate Theorem~\ref{self-dual}, we list the generators of
all self-dual ideals for ${\frak sl}_7$:\par
$\Gamma_1=\{(1,5),(2,6),(3,7)\}$,
$\Gamma_2=\{(1,4),(2,6),(4,7)\}$,
$\Gamma_3=\{(1,4),(2,5),(5,7)\}$,\par
$\Gamma_4=\{(1,3),(3,6),(4,7)\}$,
$\Gamma_5=\{(1,3),(3,5),(5,7)\}$.

\begin{rem}{Remark}   \label{nar}
The equality of Theorem~\ref{self-dual} is (almost)
an instance of the so-called ``$q=-1$
phenomenon'' studied by J.\,Stembridge \cite{stembr}.
%Consider the statistic
%\[
%   \textrm{gen}: \AD_{n} \to {\Bbb Z}_{{\ge} 0}
%\]
%that assigns to an ideal the number of its generators.
The distribution of the
statistic ``number of generators'' yields the polynomial
\[
  N_n(q)= \sum_{k=0}^{n-1} \#(\AD_{n}^{k})q^k=\sum_{k=0}^{n-1}
  \frac{1}{n}\genfrac{(}{)}{0pt}{}{n}{k}\genfrac{(}{)}{0pt}{}{n}{k+1}q^k \ ,
\]
which is often called the {\it Narayana polynomial\/}.
The $q=-1$ phenomenon is said to occur if $N_n(-1)$ counts the number of fixed
points of some natural involution on $\AD_{n}$.
We already have the involution `$\ast$' and know the number of its fixed points.
On the other hand,
it follows from \cite[Prop.\,2.2]{BSS} that
\begin{equation}   \label{-1}
N_n(-1)=\left\{ \begin{array}{rl} 0, & \textrm{ if \ $n$ \ is even} \\
(-1)^{(n-1)/2} C_{n-1/2}, &    \textrm{ if \ $n$ \ is odd} \ .
        \end{array}\right.
\end{equation}
(Actually, the authors of \cite{BSS} deal with the polynomial
$d_n(q)=(1+q)N_n(q+1)$.
However, the sign given there for the value $d_n(-2)$ should be opposite.)
Thus, we see that the $q=-1$ phenomenon  occurs up to sign.
%The meaning of sign $(-1)^{(n-1)/2}$ is however not clear.
It is interesting that Equality~\re{-1} appears also in
\cite[p.276]{reiner} in connection with a discussion of the Charney-Davis
conjecture and properties of
the Coxeter zonotope of type $A$.
\end{rem}%
The involution on $\AD_n$ (and hence on the set of Dyck paths of
semilength $n$) described in Theorem~\ref{1-1} seems to be new.

\section{{\sf ad}-nilpotent $\be$-ideals for orthogonal and
symplectic Lie algebras}    \label{classical}
\setcounter{equation}{0}

\noindent
A possible idea for constructing an involutory mapping
$\ast: \AD(\g)\to \AD(\g)$ for the other classical Lie algebras
%%$\g=\spn$ or $\son$
can be the following:

Consider the standard embedding $\g\subset {\frak sl}_N$, and choose
a Borel subalgebra $\ov{\be}\subset {\frak sl}_N$ such that
$\ov{\be}\cap\g=\be$ is a Borel subalgebra of $\g$.
Making use of the embedding $\be\subset\ov{\be}$,
one can regard $\AD(\be,\g)$ as a subset of $\AD(\ov{\be},{\frak sl}_N)$
consisting of
ideals satisfying a symmetry condition.
%(The Borel subalgebras $\be$ and $\ov{\be}$ are compatible so
%that $\ov{\be}\cap\g=\be$.)
%Therefore one can attach to an {\sf ad}-nilpotent ideal for $\g$ the
%{\sf ad}-nilpotent ideal for ${\frak sl}_N$ of the same ``shape" and
Then we apply to $\AD(\ov{\be},{\frak sl}_N)$
the duality procedure described in the previous section.
The last step should be to interpret the resulting ideal in
${\frak sl}_{N}$ as an element of $\AD(\be,\g)$.
\\[.6ex]
It turns out that this recipe yields ``expected" results for
$\spp$, but not immediately for $\sop$.
The obstacle is that the last step in the above program cannot always
be fulfilled in the orthogonal case.
%, for $\ce\subset\g\subset{\frak sl}_N$,
%the ideal $\ce^*$ may not lie to $\g$, see below.
Still, one can
modify this procedure, so that to get a suitable result for ${\frak so}_{2p+1}$.
However, I do not know how to deal with the case of ${\frak so}_{2p}$.

\begin{subs}{The symplectic case}
\end{subs}
Choose a basis for a $2p$--dimensional symplectic
$\bbk$-vector space $\VV$ so that the skew-symmetric non-degenerate bilinear form
has the matrix
 {\scriptsize
$\left(\begin{array}{cc}
 0 & \Upsilon_p \\ -\Upsilon_p & 0
 \end{array}\right)$},
where $\Upsilon_p$ is the $p\times p$ matrix whose only nonzero entries are
1's along the antidiagonal.

For any $A\in Mat_p(\bbk)$, let $\widehat A$ denote
the matrix $\Upsilon_p (A^t)\Upsilon_p$, where $A^t$ is the usual transpose
of $A$. The transformation $A \mapsto
\widehat A$ is the transpose relative to the antidiagonal.
In the above basis for $\VV$, the algebra $\spp$ has the following
block form:
\[
\spp=\{\left(\begin{array}{cr} A & B \\ C & D \end{array}
\right)\mid B={\widehat B},\ C={\widehat C},\ D=-{\widehat A} \} \ ,
\]
where $A,B,C,D$ are $p\times p$ matrices.
If $\ov{\be}$ is the standard Borel subalgebra of ${\frak sl}_{2p}$, then
$\be:=\ov{\be}\cap\spp$ is a Borel subalgebra of $\spp$.
It follows that $\AD(\spp)$ can be identified with {\it the\/} subset of
$\AD({\frak sl}_{2p})$ consisting of all Ferrers diagram that are
symmetric relative to the antidiagonal.

Let us say that $\bar I\in \AD({\frak sl}_{2p})$ is {\it self-conjugate\/},
if the corresponding Ferrers diagram is symmetric with respect
to the antidiagonal.
It is easily seen that if $\bar I\in \AD({\frak sl}_{2p})$ is self-conjugate,
then $\bar I^*$ is self-conjugate as well, see below.
This induces the desired involution on
$\AD(\spp)$, and a straightforward verification shows
that this involution satisfy properties
\re{universal} and \re{universal2}.
\\[.6ex]
Since the Ferrers diagram corresponding to an {\sf ad}-nilpotent
$\be$-ideal has a symmetry property,
we may cancel out its part which is below the antidiagonal.
What we obtain is a shifted Ferrers diagram.

\begin{figure}[htb]
\setlength{\unitlength}{0.018in}
\begin{center}
\begin{picture}(95,85)(1,0)
\multiput(10,70)(0,10){2}{\line(1,0){70}}
\put(30,60){\line(1,0){50}}
\put(50,50){\line(1,0){30}}
\multiput(60,30)(0,10){2}{\line(1,0){20}}
\multiput(70,10)(0,10){2}{\line(1,0){10}}
\qbezier[20](40,40),(62,62),(84,84)
\multiput(10,70)(10,0){8}{\line(0,1){10}}
\multiput(30,60)(10,0){6}{\line(0,1){10}}
\multiput(50,50)(10,0){4}{\line(0,1){10}}
\put(60,30){\line(0,1){20}}
\multiput(70,10)(10,0){2}{\line(0,1){40}}
%\put(70,90){\line(1,0){50}}  % \put(95,105){\vector(0,-1){9}}
%\put(90,80){\line(1,0){30}}  % \put(105,95){\vector(-1,0){9}}
%\multiput(100,40)(0,10){4}{\line(1,0){20}}
\put(-10,40){$\bar I$:}
\end{picture}
\end{center}
\caption{The self-conjugate {\sf ad}-nilpotent ideal $\bar I$ in
$\AD({\frak sl}_8)$ }
\label{sl8}
\end{figure}
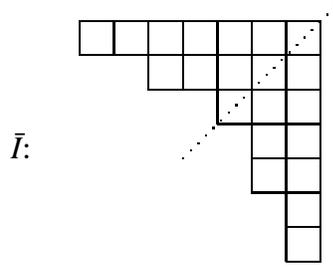

\begin{rems}{Example}  $\g={\frak sp}_8$.\\
In our matrix interpretation, the array of positive roots is
\begin{center}
$\begin{array}{ccccccc}
1000 & 1100 & 1110 & 1111 & 1121 & 1221 & 2221 \\
     & 0100 & 0110 & 0111 & 0121 & 0221        \\
     &      & 0010 & 0011 & 0021 &             \\
     &      &      & 0001 &      &      & ,
%\ap_1 & \ap_1{+}\ap_2 & \ap_1{+}\ap_2{+}\ap_3 & \ap_1{+}\ap_2{+}\ap_3{+}\ap_4 &
%\ap_1{+}\ap_2{+}2\ap_3{+}\ap_4 & \ap_1{+}2\ap_2{+}2\ap_3{+}\ap_4 &
%2\ap_1+2\ap_2+2\ap_3+\ap_4
\end{array}$
\end{center}
where the quadruple $c_1c_2c_3c_4$ stands for the root
$\sum c_i\ap_i$.
Consider the {\sf ad}-nilpotent ideal $I$ whose generators are
$\ap_1, \ap_2+\ap_3, 2\ap_3+\ap_4$. The corresponding shifted Ferrers diagram is
depicted on the left hand side in Figure~\ref{pikcha_C}.
\end{rems}
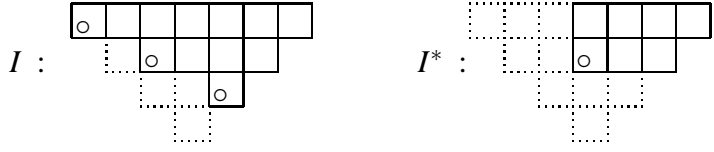
\begin{figure}[htb]
\setlength{\unitlength}{0.018in}
\begin{center}
\begin{picture}(95,45)(1,0)
\multiput(10,30)(10,0){8}{\line(0,1){10}}
\multiput(30,20)(10,0){5}{\line(0,1){10}}
\multiput(50,10)(10,0){2}{\line(0,1){10}}
\multiput(10,30)(0,10){2}{\line(1,0){70}}
\put(30,20){\line(1,0){40}}
\put(50,10){\line(1,0){10}}
%\put(70,90){\line(1,0){50}}  % \put(95,105){\vector(0,-1){9}}
%\put(90,80){\line(1,0){30}}  % \put(105,95){\vector(-1,0){9}}
%\multiput(100,40)(0,10){4}{\line(1,0){20}}
\qbezier[5](20,20),(25,20),(30,20)  % horizontal
\qbezier[10](30,10),(40,10),(50,10)     % horizontal
\qbezier[5](40,0),(45,0),(50,0)  % horizontal
\qbezier[5](20,20),(20,25),(20,30)
\qbezier[5](30,10),(30,15),(30,20)
\qbezier[10](40,0),(40,10),(40,20)
\qbezier[5](50,0),(50,5),(50,10)
\put(11,31){$\circ$}
\put(31,21){$\circ$}
\put(51,11){$\circ$}
\put(-8,20){$I\ :$}
\end{picture}
\quad\quad
\begin{picture}(95,45)(1,0)
\multiput(40,20)(10,0){4}{\line(0,1){20}}
\put(80,30){\line(0,1){10}}
%\multiput(50,10)(10,0){2}{\line(0,1){10}}
\multiput(40,30)(0,10){2}{\line(1,0){40}}
\put(40,20){\line(1,0){30}}
\qbezier[15](10,40),(25,40),(40,40)  % horizontal
\qbezier[15](10,30),(25,30),(40,30)  % horizontal
\qbezier[10](20,20),(30,20),(40,20)  % horizontal
\qbezier[15](30,10),(45,10),(60,10)     % horizontal
\qbezier[5](40,0),(45,0),(50,0)  % horizontal
\qbezier[5](10,30),(10,35),(10,40)
\qbezier[10](20,20),(20,30),(20,40)
\qbezier[15](30,10),(30,25),(30,40)
\qbezier[10](40,0),(40,10),(40,20)
\qbezier[10](50,0),(50,10),(50,20)
\qbezier[5](60,10),(60,15),(60,20)
\put(41,21){$\circ$}
\put(-5,20){$I^*\ :$}
\end{picture}
\end{center}
\caption{An {\sf ad}-nilpotent ideal in $\AD({\frak sp}_8)$ and its dual}
 \label{pikcha_C}
\end{figure}
\noindent
The dotted lines demonstrate the positive roots that are not in
$I$, and the whole array corresponds to $\Delta^+$ (or $\ut$).
The boxes marked with `$\circ$' represent the generators.
The corresponding self-conjugate ideal $\bar I \in
\AD({\frak sl}_8)$
is depicted in Figure~\ref{sl8}, where the dotted line is the antidiagonal.

\noindent
{From} the picture representing $\bar I$, we find that
$X(\bar I)=\{1,2,3,5,7\}$ and $Y(\bar I)=\{2,4,6,7,8\}$.
Therefore $X(\bar I^*)=\{2,4\}$ and $Y(\bar I^*)=\{5,7\}$.
This leads to the diagram depicted on the right hand side in
Figure~\ref{pikcha_C}. The sole generator of the ideal
$I^*$ is $\ap_2+\ap_3+\ap_4$.

Formally, our recipe for constructing the dual ad-nilpotent ideal
in $\AD(\spp)$ is as follows. We use the same coordinate system as in the
$\sln$-case. The shifted Ferrers diagram (as in Figure~\ref{pikcha_C})
is determined by the coordinates of the boxes that contain its
southwest corners, and these boxes give rise to the generators of the
respective {\sf ad}-nilpotent ideal. Suppose
$\Gamma=\{(i_1,j_1),\dots,(i_k,j_k)\}$ is the set of generators of
$I\in\AD(\spp)$, and $i_1< i_2< \ldots < i_k$. Then
$i_l<j_l$ for all $l$, $j_1< j_2< \ldots < j_k$, and $i_k+j_k\le 2p+1$.
Conversely, if a set $\Gamma$ satisfies all these inequalities, then
it is the set of generators of an {\sf ad}-nilpotent ideal.
Denoting by $\bar I$ the corresponding self-conjugate ideal in
$\AD({\frak sl}_{2p})$, we obtain

$X(\bar I)= (i_1,\dots,i_k,2p+1-j_k,\dots, 2p+1-j_1)$,

$Y(\bar I)= (j_1,\dots,j_k,2p+1-i_k,\dots, 2p+1-i_1)$.

\noindent [If $i_k+j_k=2p+1$, then one should cancel out the repetition
in the middle.]  The coordinates of vectors $X(\bar I), Y(\bar I)$ can be
paired so that the sum in each pair is equal to $2p+1$. Therefore the same
property holds for the shifted complements $X(\bar I^*), Y(\bar I^*)$.
That is, $\bar I^*$ is again a self-conjugate ideal in
$\AD({\frak sl}_{2p})$, and we can define the ideal $I^*\in \AD(\spp)$.
\\[.6ex]
Notice that $\#\Gamma(I)+\#\Gamma(I^*)=p$ and
the multiset
$\{ \Gamma(I),\Gamma(I^*)\}$ contains a unique long root, i.e.,
the distribution of long and short roots is always the same as in $\Pi$.
(A long root corresponds to the generator $(i_k,j_k)$ with
$i_k+j_k=2p+1$.) In particular, the equality $I= I^*$ is impossible, i.e.,
there are {\it no} self-dual {\sf ad}-nilpotent ideals.
\begin{rems}{Example} $\g={\frak sp}_6$. \\
In Table~\ref{tablo_C3},
we list all pairs of dual {\sf ad}-nilpotent ideals including the ideals with
one and two generators.  The column with $I$ (resp. $I^*$) contains all
ideals with one (resp.)  two generators. The numeration of simple roots is standard:
$\ap_1=\esi_1-\esi_2,\, \ap_2=\esi_2-\esi_3,\, \ap_3=2\esi_3$.
\end{rems}
\begin{table}[htb]
\begin{center}
\begin{tabular}{c|rr}
no. & $\Gamma(I)$  & $\Gamma(I^*)$ \\ \hline
1--3 & $\ap_i$ &  $\Pi\setminus\{\ap_i\}$ \\
4    & $\ap_1+\ap_2$        & $\ap_1+\ap_2, \ap_3$ \\
5    & $\ap_2+\ap_3$        & $2\ap_2+\ap_3, \ap_1$ \\
6    & $2\ap_2+\ap_3$       & $\ap_2+ \ap_3, \ap_1$ \\
7    & $\ap_1+ \ap_2+ \ap_3$ & $\ap_1+\ap_2, 2\ap_2+\ap_3$ \\
8    & $\ap_1+ 2\ap_2+\ap_3$ & $\ap_1+\ap_2+\ap_3, 2\ap_2+\ap_3$ \\
9    & $2\ap_1+2\ap_2+\ap_3$ & $\ap_1+\ap_2, \ap_2+\ap_3$
\end{tabular}
\vskip1ex
\caption{Pairs of dual {\sf ad}-nilpotent ideals in ${\frak sp}_6$}   \label{tablo_C3}
\end{center}
\end{table}
\noindent
It is clearly seen that properties of Lemmas~\ref{universal} and \ref{universal2}
are satisfied here.

\begin{subs}{The orthogonal case}
\end{subs}  \setcounter{subsubsection}{0}
Choose a basis for an $n$--dimensional orthogonal
$\bbk$-vector space $\VV$ so that the symmetric non-degenerate bilinear form
has the matrix $\Upsilon_n$.
In the above basis for $\VV$, we have:
\[
\son=\{ A \mid A=-{\widehat A} \} \ .
\]
Here we also have $\be:=\ov{\be}\cap\son$ is a Borel subalgebra.
This means that to any {\sf ad}-nilpotent $\be$-ideal in
$\son$, one can again attach a self-conjugate
{\sf ad}-nilpotent $\ov{\be}$-ideal in $\sln$.
But unlike  the symplectic case this mapping is not onto.
The reason is that the orthogonal matrices have zero antidiagonal entries.
Therefore a self-conjugate
{\sf ad}-nilpotent ideal in $\sln$ having a generator on the antidiagonal
cannot correspond to a $\be$-ideal in $\son$.
It may happen that, for $I\in\AD(\son)$, the last element
in the sequence $I \to \bar I \to \bar I^*$
cannot be interpreted as an ideal in $\son$. So, a naive attempt to
repeat the ``symplectic" procedure fails.

In the odd-dimensional case, this difficulty can be circumvented by
associating to
a $\be$-ideal in ${\frak so}_{2p+1}$  the ideal in $\spp$ having the same
shape (shifted Ferrers diagram). This is achieved by cancelling
out from a symmetric Ferrers diagram
both the antidiagonal (which corresponds to zero entries in the matrix)
and the part below the antidiagonal. This leads to a satisfactory procedure.

\begin{rems}{Example} $\g={\frak so}_7$. \\
In Table~\ref{tablo_B3},
we list all pairs of dual {\sf ad}-nilpotent ideals including the ideals with
one and two generators.  The column with $I$ (resp. $I^*$) contains all
ideals with one (resp.)  two generators. The numeration of simple roots is standard:
$\ap_1=\esi_1-\esi_2,\, \ap_2=\esi_2-\esi_3,\, \ap_3=\esi_3$. One can see
some small distinctions from Table~\ref{tablo_C3}.
\end{rems}
\begin{table}[htb]
\begin{center}
\begin{tabular}{c|rr}
no. & $\Gamma(I)$  & $\Gamma(I^*)$ \\ \hline
1--3 & $\ap_i$ &  $\Pi\setminus\{\ap_i\}$ \\
4    & $\ap_1+\ap_2$        & $\ap_1+\ap_2, \ap_3$ \\
5    & $\ap_2+\ap_3$        & $\ap_2+2\ap_3, \ap_1$ \\
6    & $\ap_2+2\ap_3$       & $\ap_2+ \ap_3, \ap_1$ \\
7    & $\ap_1+ \ap_2+ \ap_3$ & $\ap_1+\ap_2, \ap_2+2\ap_3$ \\
8    & $\ap_1+ \ap_2+2\ap_3$ & $\ap_1+\ap_2+\ap_3, \ap_2+2\ap_3$ \\
9    & $\ap_1+2\ap_2+2\ap_3$ & $\ap_1+\ap_2, \ap_2+\ap_3$
\end{tabular}
\vskip1ex
\caption{Pairs of dual {\sf ad}-nilpotent ideals in ${\frak so}_7$}   \label{tablo_B3}
\end{center}
\end{table}
\noindent
Again, the properties of Lemmas~\ref{universal} and \ref{universal2}
are satisfied here.
In the following section, we also summarize some other properties of the
duality mapping that
are inspired by our computations in classical cases.

\section{Towards the general case}    \label{vmeste}
\setcounter{equation}{0}

\noindent
In view of Theorem~\ref{1-1}, %and Lemma~\ref{universal},
it is natural to ask whether there is
a natural involutory mapping
 $\ast : \AD(\g) \to \AD(\g)$ for any simple
Lie algebra $\g$ such that \\
\centerline{$\#(\Gamma(I))+\#(\Gamma(I^*))=\rk\g$}
and the two properties of
Lemmas~\ref{universal} and \ref{universal2} are also satisfied\,?
\par
It is plausible that a conjectural definition of duality should exploit
somehow admissible elements of $\widehat W$ and the simplex $\tilde D$.
Although my attempts to define such a mapping in a uniform way were unsuccessful, I believe that such a mapping does exist.
\par
%%This section consists essentially of various speculations related to
%%a conjectural duality on the set of {\sf ad}-nilpotent ideals.
Since an {\sf ad}-nilpotent ideal $I\in\AD(\g)$ is completely determined
by the corresponding antichain $\Gamma=\Gamma(I)\subset\Delta^+$,
properties of the conjectural duality
on $\AD(\g)$ can be restated in terms of antichains in $\Delta^+$.
Let $\AN(\Delta^+)$ denote the set of all antichains in  $\Delta^+$.
For a moment, we assume that $\Delta$ is not necessarily irreducible,
and $\Delta=\underset{i}{\sqcup} \Delta_i$, where each $\Delta_i$ is an
irreducible root system and the rank of $\Delta_i$ is $p_i$.
\begin{s}{Conjecture}  \label{main-conj}
There exists a natural involutory mapping
\[
   \ast:  \AN(\Delta^+) \to \AN(\Delta^+)
\]
such that the following holds for\/ $\Gamma\in \AN(\Delta^+)$:
\begin{itemize}
\item[\sf (i)]  $\Gamma^*=\sqcup (\Gamma\cap\Delta_i)^*$ and
$(\Gamma\cap\Delta_i)^*$ depends only on $\Gamma\cap\Delta_i$;
\item[\sf (ii)]
$\#(\Gamma\cap\Delta_i) +\#(\Gamma^*\cap\Delta_i)=p_i$ for all $i$;
\item[\sf (iii)]   Suppose $\Gamma$ contains a simple root $\ap$. Write
$\Delta(\Pi\setminus\{\ap\})$ for the root subsystem
spanned by the set of simple roots $ \Pi\setminus\{\ap\}$.
Then $\Gamma^*\subset \Delta(\Pi\setminus\{\ap\})^+$ and
moreover, $\Gamma^*=(\Gamma\setminus\{\ap\})^*$, where
$\Gamma\setminus\{\ap\}$ is regarded as antichain
in $\Delta(\Pi\setminus\{\ap\})^+$;
\item[\sf (iv)]  ($\approx$ a converse to the previous property)
\ If $\Gamma\subset \Delta(\Pi\setminus\{\ap\})^+$, then \\
$\Gamma^*=\{\ap\}\cup \{
\text{the dual of $\Gamma$ taken in }
 \Delta(\Pi\setminus\{\ap\})^+\}$;
\item[\sf (v)]   If $\Delta$ is irreducible, then $(\Delta^+(k))^*=\Delta^+(h+1-k)$,
where $h$ is the Coxeter number of $\Delta$
(cf. Lemma~\ref{universal2}).
\item[\sf (vi)]  the distribution of long and short roots in the multiset
$\{\Gamma,\Gamma^*\}$ is the same as in $\Pi$.
(This condition is vacuous in the simply-laced case)
\end{itemize}
\end{s}%
It is easy to see that the duality defined for the root systems of type
$\GR{A}{p},\,\GR{B}{p},\,\GR{C}{p}$ satisfies
all these properties. Also, it is immediate that `$\ast$' can uniquely be
defined for $\GR{G}{2}$.

\noindent
Now, we again assume that $\Delta$ is irreducible.
Clearly, a necessary condition for such a duality to exist is that
the number of antichains
of cardinality $k$ ought to be equal to the number of antichains of
cardinality $p-k$.
%%(Recall that $p=\rk\g$.)
This holds in all cases, where the corresponding values are known,
see below.
If $k=0$, then the assertion follows from
Proposition~\ref{vershiny}. In case  $k=1$, one should be able to
prove that the number of
positive roots is equal to the number of antichains of
cardinality $p-1$. Unfortunately, the only proof I know amounts to
a case-by-case verification,
%even this simply-formulated question
%has no general answer as yet.

For each simple Lie algebra $\g$, we define an analogue
of Narayana polynomial as follows. Let $d_k(\g)$ be the
number of all {\sf ad}-nilpotent ideals with $k$ generators
or, equivalently, the number of all $k$-element antichains
in $\Delta^+$.
Then
\begin{equation}  \label{gen-nar}
   \N_{\g}(q)=\sum_{i=0}^{p} d_k(\g) q^k
\end{equation}
is said to be the {\it Narayana polynomial of type} $\g$ (or,
a generalized Narayana polynomial). Clearly, $d_0(\g)=d_p(\g)=1$ and
$d_1(\g)=\#\Delta^+$. By Theorem~\ref{main2}, $d_{p-1}(\g)$
%%the coefficient of $q^{p-1}$
equals the number of integral points lying on the edges of the simplex
$\tilde D$ (except of the unique integral vertex).
Below, we list
all known to us generalized Narayana polynomials:
\begin{itemize}
\item[]   $\N_{A_p}(q)=\displaystyle
\sum_{k=0}^p \frac{1}{p+1}
\genfrac{(}{)}{0pt}{}{p+1}{k}\genfrac{(}{)}{0pt}{}{p+1}{k+1}q^k$;
\item[]  $\N_{B_p}(q)=\N_{C_p}(q)=\displaystyle
\sum_{k=0}^p \genfrac{(}{)}{0pt}{}{p}{k}^2 q^k$;
\item[]   $\N_{G_2}(q)=1 +  6q+      q^2$;
\item[]   $\N_{F_4}(q)=1 +  24q  + 55q^2 + 24 q^3 +q^4$ ;
\item[]   $\N_{E_6}(q)=1 +  36q  +204q^2 + 351q^3 + 204^4 +   36q^5 +    q^6$ ;
\item[]   $\N_{E_7}(q)=1 +  63q  +546q^2 +1470q^3 +1470^4 +  546q^5 +  63q^6+   q^7$ ;
\item[]   $\N_{E_8}(q)=1 + 120q +1540q^2 +6120q^3 +9518^4 + 6120q^5 +1540q^6+120q^7+ q^8$ .
\end{itemize}
In type $A$, it is the usual Narayana polynomial (cf. \ref{nar}).
The result for types $B$ and $C$ follows from \cite[Corollary\,5.8]{grek}.
In that place, Athanasiadis computes the number of non-nesting partitions on
$\GR{B}{p}$ or $\GR{C}{p}$ whose `type' has $k$ parts.
However, it follows from his previous exposition that a non-nesting partition
whose type has $k$ parts is exactly an antichain of cardinality $p-k$.
%%In exceptional cases we performed explicit calculations.
The case of $\GR{G}{2}$ is trivial and that of $\GR{F}{4}$ is relatively
easy.
\\
The case of $\GR{E}{n}$ requires more work.
The result can be obtained through the counting all integral points in $\tilde D$ and use of
Theorem~\ref{main2}.
%To compute the coefficient of $q^3$, we use some ideas of Shi,
%see \cite[3.4]{shi}, and for the coefficient of $q^5$, we find
%all integral points on the edges of $\tilde D$.
\\[.5ex]
I also conjecture that
$\N_{D_p}(q)=\displaystyle
\sum_{k=0}^p \Bigl(\genfrac{(}{)}{0pt}{}{p}{k}^2-
\frac{p}{p-1}\genfrac{(}{)}{0pt}{}{p-1}{k}\genfrac{(}{)}{0pt}{}{p-1}{k-1}
\Bigr) q^k$.  \\
%% and  \\
%%$N_{E_6}(q)=1+36 q + 203 q^2+353 q^3+ 203 q^4+ 36q^5+ q^6$.
In fact, $\genfrac{(}{)}{0pt}{}{p}{k}^2-
\frac{p}{p-1}\genfrac{(}{)}{0pt}{}{p-1}{k}\genfrac{(}{)}{0pt}{}{p-1}{k-1}$
is the number of {\sl non-crossing\/}
partitions on $\GR{D}{p}$ whose type has $k$ parts \cite[Sect.\,4]{victor}.
And I hope that the similarity between the non-crossing and non-nesting
partitions known for $\GR{A}{p},\GR{B}{p}$ and $\GR{C}{p}$ remains true
also for $\GR{D}{p}$.
Thus, all known generalized Narayana polynomials are palindromic.

By \cite{CP2},
we have $\N_\g(1)=\#\AD(\g)=\displaystyle
\prod_{i=1}^p \frac{h+e_i+1}{e_i+1}$.
%%, where the $e_i$'s are the exponents and $h$ is the Coxeter number
%%of $\g$. It is also the value.
It would be interesting to find a uniform expression
for the coefficients of the generalized Narayana polynomials.

Another intriguing feature is that there are nice formulae for the values
$\N_\g(-1)$. For $\GR{A}{p}$, we refer again to \ref{nar}.
The $\GR{B}{p}$- or $\GR{C}{p}$-case amounts to a well-known combinatorial identity:
\[
    \sum_{k=0}^p (-1)^k \genfrac{(}{)}{0pt}{}{p}{k}^2=
   \left\{ \begin{array}{cl}  0, & \textrm{ if $p$ is odd}, \\
              (-1)^{p/2}\genfrac{(}{)}{0pt}{}{p}{p/2},  & \textrm{ if $p$ is even}.
            \end{array} \right.
\]
Combining the expressions for $\GR{A}{p}$ and $\GR{B}{p}$ cases, we conclude
that if the conjectural formula for $\N_{D_p}(q)$ is true, then
\[
    \N_{D_p}(-1)=\left\{ \begin{array}{cl}  0, & \textrm{ if $p$ is odd}, \\
              (-1)^{p/2}\Bigl[ \genfrac{(}{)}{0pt}{}{p}{p/2}-
       2\genfrac{(}{)}{0pt}{}{p-2}{p/2-1}\Bigr]=
        (-1)^{p/2}2\genfrac{(}{)}{0pt}{}{p-2}{p/2},
              & \textrm{ if $p$ is even}.
            \end{array} \right.
\]
One may also observe that if $p$ is even, then
$(-1)^{p/2} \N_{\g}(-1)$ is positive for all examples.

\end{document}